\documentclass[aos,MSNbibl,nameyear,dvips]{arximspdf}
\usepackage{dcolumn}
\usepackage{graphicx}

\doi{10.1214/11-AOS932}
\volume{39}
\issue{6}
\pubyear{2011}
\firstpage{3182}
\lastpage{3210}

\makeatletter

\newcolumntype{d}[1]{D{.}{.}{#1}}

\newtheorem{proposition}{Proposition}
\newtheorem{corollary}{Corollary}

\newproclaim{definition}{Definition}

\newtheorem{theorem}{Theorem}

\newcommand{\dist}{\operatorname{dist}}
\renewcommand{\vec}{\operatorname{vec}}

\newcommand{\indep}{\perp\!\!\!\!\perp}

\newcommand{\var}{\operatorname{var}}
\newcommand{\rank}{\operatorname{rank}}
\newcommand{\diag}{\operatorname{diag}}
\newcommand{\spc}{{\mathcal S}}
\newcommand{\bblfnull}{{\mathbf0}}

\newcommand{\cov}{\operatorname{cov}}
\renewcommand{\sp}{\operatorname{span}}

\newcommand{\cbf}{{\mathbf c}}
\newcommand{\ebf}{{\mathbf e}}
\newcommand{\sbf}{{\mathbf s}}
\newcommand{\wbf}{{\mathbf w}}
\newcommand{\ybf}{{\mathbf y}}
\newcommand{\Abf}{{\mathbf A}}
\newcommand{\Bbf}{{\mathbf B}}
\newcommand{\Dbf}{{\mathbf D}}
\newcommand{\Fbf}{{\mathbf F}}
\newcommand{\Hbf}{{\mathbf H}}
\newcommand{\Ibf}{{\mathbf I}}
\newcommand{\Jbf}{{\mathbf J}}
\newcommand{\Kbf}{{\mathbf K}}
\newcommand{\Mbf}{{\mathbf M}}
\newcommand{\Pbf}{{\mathbf P}}
\newcommand{\Qbf}{{\mathbf Q}}
\newcommand{\Tbf}{{\mathbf T}}
\newcommand{\Ubf}{{\mathbf U}}
\newcommand{\Vbf}{{\mathbf V}}
\newcommand{\Xbf}{{\mathbf X}}
\newcommand{\Zbf}{{\mathbf Z}}
\newcommand{\alphabf}{{\bolds\alpha}}
\newcommand{\betabf}{{\bolds\beta}}
\newcommand{\deltabf}{{\bolds\delta}}
\newcommand{\zetabf}{{\bolds\zeta}}
\newcommand{\etabf}{{\bolds\eta}}
\newcommand{\thetabf}{{\bolds\theta}}
\newcommand{\xibf}{{\bolds\xi}}
\newcommand{\phibf}{{\bolds\phi}}
\newcommand{\Lambdabf}{{\bolds\Lambda}}
\newcommand{\Sigmabf}{{\bolds\Sigma}}
\newcommand{\Upsilonbf}{{\bolds\Upsilon}}
\newcommand{\Psibf}{{\bolds\Psi}}
\newcommand{\vb}{{\mathbf v}}

\newcommand{\R}{\mathbb R}
\newcommand{\Xb}{\mathbf X}
\newcommand{\trans}{^{\top}}

\newcommand{\inv}{^{-1}}
\newcommand{\cid}{\stackrel{{\mathcal D}}{\longrightarrow}}

\newcommand{\hx}{{\H}}
\newcommand{\hxstar}{{\H'}}
\newcommand{\half}{^{1/2}}
\renewcommand{\H}{{\mathcal H}}
\newcommand{\zero}{\mathbf{0}}
\newcommand{\one}{\mathbf{1}}
\newcommand{\Psub}{{P}}
\newcommand{\mhalf}{^{-{1/2}}}
\newcommand{\cip}{\stackrel{\Psub}{\to}}

\makeatother

\begin{document}
\begin{frontmatter}

\title{Principal support vector machines for linear and nonlinear
sufficient dimension reduction}
\runtitle{Hyperplane alignment}

\begin{aug}
\author[A]{\fnms{Bing} \snm{Li}\corref{}\thanksref{t1}\ead[label=e1]{bing@stat.psu.edu}},
\author[B]{\fnms{Andreas} \snm{Artemiou}\ead[label=e2]{aartemio@mtu.edu}}
\and
\author[C]{\fnms{Lexin} \snm{Li}\thanksref{t2}\ead[label=e3]{li@stat.ncsu.edu}}
\runauthor{B. Li, A. Artemiou and L. Li}
\affiliation{Pennsylvania State University, Michigan Technological
University and~North Carolina State University}
\address[A]{B. Li\\
Department of Statistics\\
Pennsylvania State University\\
326 Thomas Building\\
University Park, Pennsylvania 16802\\
USA\\
\printead{e1}} 
\address[B]{A. Artemiou\\
Department of Mathematical Sciences\\
Michigan Technological University\\
Fisher Hall, Room 306\\
1400 Townsend Drive\\
Houghton, Michigan 49931\\
USA\\
\printead{e2}}
\address[C]{L. Li\\
Department of Statistics\\
North Carolina State University\\
2311 Stinson Drive\\
Campus Box 8203\\
Raleigh, North Carolina 27695-8203\\
USA\\
\printead{e3}}
\end{aug}

\thankstext{t1}{Supported in part by
NSF Grants DMS-07-04621 and DMS-08-06058.}

\thankstext{t2}{Supported in part by NSF Grant DMS-11-06668.}

\received{\smonth{10} \syear{2010}}
\revised{\smonth{9} \syear{2011}}

%
\begin{abstract}
We introduce a principal support vector machine (PSVM) approach that
can be used for both linear and nonlinear sufficient dimension
reduction. The basic idea is to divide the response variables into
slices and use a modified form of support vector machine to find the
optimal hyperplanes that separate them. These optimal hyperplanes are
then aligned by the principal components of their normal vectors. It is
proved that the aligned normal vectors provide an unbiased, $\sqrt
n$-consistent, and asymptotically normal estimator of the sufficient
dimension reduction space. The method is then generalized to nonlinear
sufficient dimension reduction using the reproducing kernel Hilbert
space. In that context, the aligned normal vectors become functions and
it is proved that they are unbiased in the sense that they are
functions of the true nonlinear sufficient predictors. We compare PSVM
with other sufficient dimension reduction methods by simulation and in
real data analysis, and through both comparisons firmly establish its
practical advantages.
\end{abstract}

%

\end{frontmatter}

\section{Introduction}
With the increase of computer power in storing and processing data,
high dimensional data have become increasingly prevalent across
many disciplines. The demand for effective methods to extract
useful information from such data has led inevitably to dimension
reduction, an area that has undergone tremendous development
during the past two decades.

Let $\Xbf$ be a $p$-dimensional predictor and $Y$ be a response variable.
In its classical form, sufficient dimension reduction (SDR)
[Li (\citeyear{Li91N2}, \citeyear{Li92}), \citet{CoWe91},
\citet{Coo98}] identifies a small number of
linear combinations of predictors that can replace the original predictor
vector~$\Xb$ without loss of information on the conditional distribution
of $Y$ given $\Xb$. In other words, the objective is to find a
$p \times d$ ($d < p$) matrix $\etabf$ such that the following
conditional independence holds:
%
\begin{equation}\label{relmodel}
Y \indep\Xb| \etabf\trans\Xb.
\end{equation}
In this relation, the identifiable parameter is the subspace spanned by
the columns of $\etabf$, rather than $\etabf$ itself. The intersection
of all subspaces satisfying~(\ref{relmodel}), provided itself satisfies
(\ref{relmodel}), is called the central subspace, and is denoted by
$\spc_{Y|\Xbf}$ [\citet{Coo94}]. \citet{Coo96} and
\citet{YinLiCoo08} showed that $\spc_{Y|\Xbf}$ uniquely exists
under very mild conditions. Thus, we assume its existence throughout
this article. Many methods have been proposed for this problem since
the publication of the original works. See, for example,
\citet{CooLi02}, \citet{Xiaetal02}, \citet{YinCoo02},
\citet{Funetal02}, \citet{LiZhaChi05},
\citet{CooNi05}, \citet{LiWan07}, \citet{LiDon09}.

A more general sufficient dimension reduction problem, as formulated in
\citet{Coo07}, is to seek an arbitrary function $\phibf\dvtx\R^p \to\R^d$
such that
%
\begin{equation}\label{eqnonlinear}
Y \indep\Xbf| \phibf(\Xbf).
\end{equation}
We refer to this problem as \textit{nonlinear sufficient dimension reduction},
and any one-to-one function of $\phibf(\Xbf)$ as the \textit{nonlinear
sufficient predictor}.
Several
recent works pioneered estimation procedures for nonlinear dimension
reduction of this type, including \citet{Wu08}, \citet{WuLiaMuk},
\citet{Wan08} and \citet{YehHuaLee09}, by extending sliced inverse
regression [SIR; \citet{Li91N2}] from different angles.

In this paper, we propose a sufficient dimension reduction method, to be
called the principal support vector machine (PSVM), that can extract
the sufficient
predictors in both problems (\ref{relmodel}) and (\ref{eqnonlinear}).
Let $(\Xbf_1,Y_1), \ldots, (\Xbf_n, Y_n)$ be a sample of $(\Xbf, Y)$.
The basic idea of PSVM is to divide $\Xbf_1, \ldots,
\Xbf_n$ into several slices according to the values of the responses,
and then use support vector machine
[SVM; \citet{Vap98}] to
find the optimal hyperplanes that separate these slices. The optimal
hyperplanes are then aligned by applying principal component analysis
to their normal vectors. We show that the principal components are,
in fact, an unbiased estimator of the central subspace~$\spc_{Y|\Xbf}$.
This idea is then extended to the nonlinear dimension reduction problem
(\ref{eqnonlinear}) via the reproducing kernel Hilbert space [RKHS;
\citet{Aro50}, \citet{HsiRen09}]. In this context, the normal vectors
in the linear case become functions in the RKHS. It is shown that the normal
functions thus derived are functions of $\phibf$ in the general problem
(\ref{eqnonlinear}). This is, to our knowledge, the first result of this
type.

Our proposal is noticeably different from the existing SDR methods in the
following respects. First, PSVM is developed under, and for,\vadjust{\goodbreak} a unified
framework of linear and nonlinear sufficient dimension reduction. Such
a~standpoint allows us to formulate some theoretical properties, such as
unbiasedness, more rigorously and generally than previous works. Second,
PSVM improves the accuracy for sufficient dimension reduction, for the
following reason. It is well known that a regression surface is more
accurately estimated at the center of the data cloud than at the outskirt.
However, an inverse regression based method, such as SIR, tends to
downweight the slice means near the center due to their shorter lengths.
Since PSVM relies on separating hyperplanes rather than slice means, it
makes better use of the central portion of the data than inverse
regression. This improvement is clearly demonstrated in our numerical
studies. Finally, PSVM establishes a~firm connection between sufficient
dimension reduction and the acclaimed machine learning technique, support
vector machine, both of which have been extensively used in high dimensional
data analysis. This combination brings fresh insights and further advances
to both subjects. Along with the theoretical development of PSVM, we develop
a more complete asymptotic theory for SVM than previously given, and
introduce the notion of invariant kernel for SVM. Meanwhile, we expect
some inherent advantages of SVM to benefit sufficient dimension
reduction estimation.
For instance, SVM tends to be more robust against outliers than a
typical moment method.
This is because the separating hyperplanes are largely determined by
the support vectors lying in the interior of the data cloud, as a result
an observation far away from the data cloud has
less influence than a typical moment-based estimator. In this sense,
SVM behaves more like a median than a~mean.
It is also expected to help address several challenging issues facing the
existing SDR methods, such as small-$n$--large-$p$ and presence of
categorical predictors. However, due to limited space these potential
advantages cannot be fully discussed within this paper. Some of them,
such as robustness and categorical predictors, are further explored in
\citet{Art10}.

The rest of the paper is organized as follows. In Section
\ref{sectionhyperplanealignment}, we illustrate the basic idea of
PSVM by examples and figures, and give intuitions about why it works.
In Section \ref{sectionlinearha}, we formally
introduce the linear PSVM and study its population-level properties
in terms of its unbiasedness as an estimator of the central subspace.
In Section \ref{sectionestimationprocedure}, we develop the estimation
procedures for the linear PSVM, and describe how to implement it using
standard SVM packages. In Section \ref{sectionnonlinearha}, we
generalize the linear PSVM to the kernel PSVM to solve the nonlinear
sufficient dimension reduction problem, and establish its unbiasedness
in this general setting. In Section~\ref{sectionestimationkernelha},
we develop an algorithm to implement the kernel PSVM, and introduce the
notion of invariant kernel. In Section \ref{sectionasymptoticanalysis}, we study the asymptotic properties for the linear PSVM
estimator.
Though the identified subspaces are asymptotically
consistent, they are almost surely incorrect for finite sample sizes.
Thus, in Section \ref{sectionsimulationstudies}, we compare the
linear and kernel PSVM with other dimension reduction methods in finite
sample by simulation.
In Section \ref{sectiondataanalysis}, we apply it to analyze a data set
concerning the recognition of vowels, and make further comparisons in
the practical setting. All the proofs are given in a~complementary document
published online by \textit{The Annals of Statistics}.

\section{Principal support vector machine: The basic idea}
\label{sectionhyperplanealignment}

The idea of the principal support vector machine arises from an
interplay of several ideas:
sliced inverse regression, support vector machine, and contour regression
[\citet{LiZhaChi05}]. In this section, we illustrate this idea
by two simple examples that cover both linear and nonlinear dimension
reduction. Throughout this paper, $\Xbf$ represents a random vector; $X_r$
represents the $r$th component of $\Xbf$; $\Xbf_i$ represents the $i$th
random vector from a sample~$\Xbf_1, \ldots, \Xbf_n$, and $X_{ir}$
represents the $r$th component of $\Xbf_i$.

First, consider the regression model
%
\begin{equation}\label{eqsimplelinearmodel}
Y = f(X_1 + 2 X_2) + \varepsilon,
\end{equation}
where $\varepsilon\indep(X_1, X_2)$.
This is a linear sufficient dimension reduction problem, in which the central
subspace is spanned by $(1, 2)\trans\in\R^2$. Note that the contours for
the regression function is the set $\{(x_1, x_2)\dvtx x_1 + 2 x_2 = c \}$,
which is uniquely associated with the vector $(1,2)\trans$. Based on this
intuition, \citet{LiZhaChi05} introduced the contour regression,
which estimates the contour directions by the directions in $\Xbf$
that are
aligned with the smallest increments in $Y$.

Here, we propose to identify the contours by the separating
hyperplanes derived from the support vector machine as applied to
different slices of~$\Xbf$, formed according to the values of $Y$.
Let
$
S_1 = \{\Xbf_i\dvtx Y_i \le c \}$ and $S_2 = \{ \Xbf_i \dvtx\allowbreak Y_i > c \}
$
for some constant $c$.
We use SVM to obtain the optimal separating
hyperplane of $S_1$ and~$S_2$, and repeat the process
to obtain several hyperplanes. Intuitively, the normals of these
hyperplanes are roughly aligned with the directions in which the
regression surface varies---directions that form
the central subspace. We use the principal components of these
normals to estimate the central subspace. A related idea is \citet{Loh02}, who
proposed to divide each individual predictor according to the mean of $Y$
and assess the importance of that predictor by its degree of separation.

\begin{figure}

\includegraphics{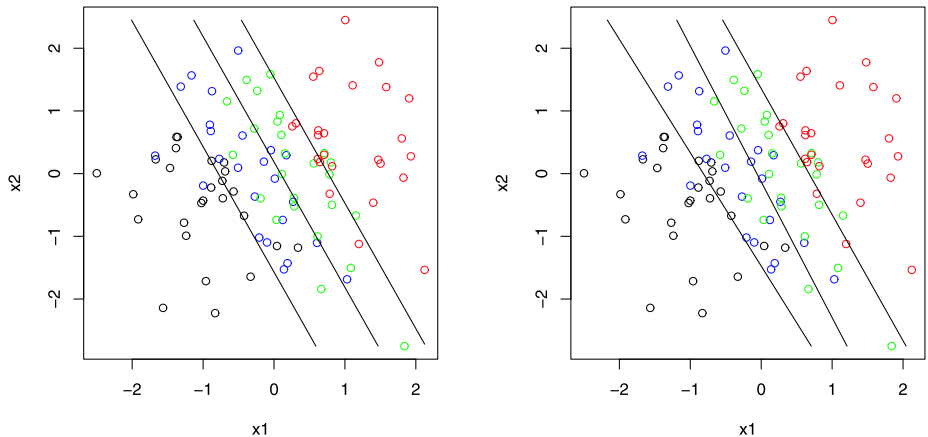}

\caption{Linear contours for model $Y = 2 X_1 + X_2 + \varepsilon$.
Left panel: true contours; right panel: contours based on linear SVM.
Contour levels are three evenly spaced sample quantiles of
$Y_1, \ldots, Y_n$.
The sample size is $n = 100$.}
\label{figurelinearcontour}
\end{figure}

As an illustration, we generate 100 replications from model
(\ref{eqsimplelinearmodel}) where~$f$ is
taken to be the identity mapping. We divide $\Xbf_1, \ldots,
\Xbf_{100}$ into 4 slices according to the 25th, 50th, 75th
sample quantiles of $Y_1, \ldots, Y_n$, as
indicated in Figure \ref{figurelinearcontour} by differently colored
dots. Application of SVM between these
slices yields three hyperplanes, represented by the solid
lines on the right panel, which closely resemble the contours
derived from the true model, as shown on the left panel.
Clearly, the normals of the three hyperplanes give close
estimate of the central subspace.

We can apply the same idea to sufficient nonlinear dimension
reduction. Let
%
\begin{equation}\label{eqnonlinearexample}
Y = f(X_1 + X_2^2) + \varepsilon,
\end{equation}
where $f$ is an unknown function. The contours of this function
are of the form $\{(x_1, x_2)\dvtx x_1 + x_2^2 = c \}$, which are
no longer hyperplanes in $\R^2$.
However, if we map ${\mathbf x}$ to a higher dimensional space of
functions of ${\mathbf x}$ that is rich enough to contain $x_1 +
x_2^2$, then
the contours become hyperplanes
again. We can apply SVM at that level to find the optimal
hyperplanes, and then map them back to the ${\mathbf x}$-space to
extract the nonlinear predictor.
Usually, in conjunction with mapping a low-dimensional
regressor to a high-dimensional regressor, a Tikhonov-type
regularization is applied, so
that the overfitting tendency of increased dimension is counteracted by
the regularization.

\begin{figure}

\includegraphics{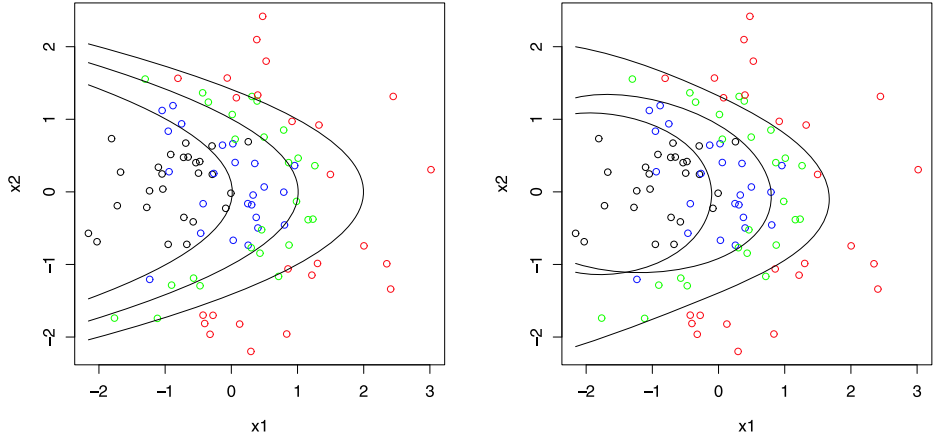}

\caption{Nonlinear contours for model $Y = X_1 + X_2^2 + \varepsilon$.
Left panel: true contours; right panel: contours based on kernel SVM
with gauss radial kernel.
Contour levels are three evenly spaced sample quantiles of $Y_1,
\ldots, Y_n$.
The sample size is $n = 100$.}
\label{figurenonlinearcontour}
\end{figure}

As in the linear case we generate 100 replications from model (\ref
{eqnonlinearexample}) and use the same set of quantiles to slice the
response. The curves in the left panel in Figure
\ref{figurenonlinearcontour} are the true contours computed from the
function $y = x_1 + x_2^2$. Those in the right panel are obtained by
first applying kernel SVM (with Gaussian radial basis) to find
hyperplanes in $\R^{100}$ and then mapping them back to $\R^2$.
Clearly, any function of $(x_1, x_2)$ that generates the contours in
the right panel would closely resemble the true predictor $x_1 +
x_2^2$, modulo a~monotone transformation.

\section{PSVM for linear sufficient dimension reduction}
\label{sectionlinearha}

We first develop PSVM for linear sufficient dimension
reduction. We begin with a
population-level formulation of SVM, since it is usually
described at the sample level, which is not the best way to
set up our problem.
For now, assume $Y$ to be a binary random variable
taking values $-1$ and $1$. The soft-margin SVM is defined
through the following optimization:
%
\begin{eqnarray}\label{eqlinearsvmobjectivefunction}
&&\mbox{minimize }  {{\bolds\psi}}\trans{{\bolds\psi
}}+ \frac{\lambda}{n} \sum
_{i=1}^n \xi_i  \qquad\mbox{among $({{\bolds\psi}}, t, \xibf)
\in\R^p
\times\R\times\R^n$} \nonumber\\[-8pt]\\[-8pt]
&&\qquad\mbox{subject to } \xi_i \ge0,  Y_i [ {{\bolds\psi
}}\trans(\Xbf_i -
\bar\Xbf) - t ] \ge1 - \xi_i, \qquad i = 1, \ldots, n,
\nonumber
\end{eqnarray}
where $\lambda$ is a positive constant often referred to
as the ``cost.'' See Vapnik [(\citeyear{Vap98}), page 411] for the
intuitions behind this construction. If $({{\bolds\psi}}^*, t^*,
\xibf^*)$ is the solution to
(\ref{eqlinearsvmobjectivefunction}), then the
set $\{{\mathbf x}\dvtx {{{\bolds\psi}}^{*\top}}{\mathbf x}= t^*\}$
is the optimal hyperplane that separates
$\{\Xbf_i\dvtx Y_i = -1\}$ and $\{\Xbf_i\dvtx Y_i = 1 \}$.

Although this representation defines the algorithm, it does
not tell us what objective function is minimized at
the population level. To see things more clearly, we first
carry out the optimization for a fixed $({{\bolds\psi}}, t)$.
This amounts to minimizing $\sum_{i=1}^n \xi_i$ subject
to $\xi_i \ge\max\{ 0, 1 - Y_i
[ {{\bolds\psi}}\trans(\Xbf_i - \bar\Xbf) - t ] \}$. The optimal
solution is
$
\xi_i^* = \{1 - Y_i [ {{\bolds\psi}}\trans( \Xbf_i - \bar
\Xbf) - t ] \}^+
$
where $a^+ =\max(a,0)$.
Substituting $\xi_i^*$ into
(\ref{eqlinearsvmobjectivefunction}), we have
%
\begin{equation}\label{eqclue}
{{\bolds\psi}}\trans{{\bolds\psi}}+ \frac{\lambda}{n}
\sum_{i = 1}^n \{1 - Y_i [ {{\bolds\psi}}\trans(\Xbf_i -
\bar\Xbf) - t ]\}^+.
\end{equation}
This corresponds to the following objective function at the
population level:
%
\begin{equation}\label{eqclassical}
{{\bolds\psi}}\trans{{\bolds\psi}}+ \lambda E
\bigl[1 - Y \bigl( {{\bolds\psi}}\trans(\Xbf- E \Xbf) - t \bigr)\bigr]^+.
\end{equation}
The hyperplane that minimizes this criterion can be viewed
as that which best separates the \textit{conditional
distributions}\vadjust{\goodbreak}
of $\Xbf| Y = -1$ and $\Xbf| Y = 1$. \citet{JiaZhaCai08} used a slight variation of representation
(\ref{eqclassical}) to derive the asymptotic distribution of
SVM. We also use a representation similar to (\ref{eqclassical}) but
with two important
modifications, as we describe below.

Return now to sufficient dimension reduction problem
(\ref{relmodel}) where $Y$ is an arbitrary random variable (in
particular, it
can be either continuous or categorical).
Let $\Omega_Y$ be the support of $Y$ and let $A_1$ and $A_2$
be disjoint subsets of $\Omega_Y$. Let $\tilde Y$ be the
discrete random variable defined by
%
\begin{equation}\label{eqtildey}
\tilde Y = I(Y \in A_2) - I(Y \in A_1).
\end{equation}
We introduce the following objective function for linear
SDR:
%
\begin{equation}\label{eqmodified}
L({{\bolds\psi}}, t) = {{\bolds\psi}}\trans\Sigmabf
{{\bolds\psi}}+ \lambda E
\{1 - \tilde Y [ {{\bolds\psi}}\trans( \Xbf- E \Xbf) - t ] \}^+,
\end{equation}
where $\Sigmabf= \var(\Xbf)$.
Compared with (\ref{eqclassical}), we have made two modifications.
First, we allow $\tilde Y$ to take the value 0, so that we can
use a pair of disjoint subsets that are not a partition
of $\Omega_Y$. Second, we have inserted
$\Sigmabf$ in the first term of (\ref{eqmodified}). This is so
that the objective function transforms in a desired manner. We will
return to this point in Section \ref{sectionestimationkernelha}.

We now establish the unbiasedness of the normal vector for the
optimal separating hyperplane in SVM as an estimator of the
central subspace.
Let~$F_n$ be the empirical distribution based on the sample
$(\Xbf_1, Y_1), \ldots,\allowbreak (\Xbf_n, Y_n)$, $F_0$ be the true
distribution of $(\Xbf, Y)$, and $\Tbf$ be a statistic that
can be expressed as a matrix-valued function of the distribution of
$(\Xbf, Y)$.
In our context, we say $\Tbf(F_n)$ is an unbiased estimator
of $\spc_{Y|\Xbf}$, if it satisfies
%
\begin{equation}\label{eqlinearunbiased}
\sp[ \Tbf(F_0)] \subseteq\spc_{Y|\Xbf}.
\end{equation}

\begin{theorem}\label{theoremmain} Suppose
$E(\Xbf| \etabf\trans\Xbf)$ is a linear
function of $\etabf\trans\Xbf$, where $\etabf$ is as defined
in (\ref{relmodel}). If $({{\bolds\psi}}^*, t^*)$ minimizes the
objective function (\ref{eqmodified}) among all $({{\bolds\psi
}}, t)
\in\R^p \times\R$, then $ {{\bolds\psi}}^* \in\spc
_{Y|\Xbf}$.
\end{theorem}

The linearity condition on $E(\Xbf| \etabf\trans\Xbf)$ in the
theorem is well known and generally assumed in the SDR
literature. See, for example, \citet{LiDua89}, \citet{Li91N2}
and \citet{LiDon09}. It implies
%
\begin{equation}\label{eqlinearimplication}
E ( {{\bolds\psi}}\trans\Xbf| \etabf\trans\Xbf) =
{{\bolds\psi}}\trans\Pbf_{\etabf
}\trans(\Sigmabf) \Xbf,
\end{equation}
where $\Pbf_\etabf(\Sigmabf)$ is the projection matrix
$\etabf(\etabf\trans\Sigmabf\etabf)^{-1} \etabf\trans\Sigmabf$
[\citet{Coo98}]. It is satisfied when $\Xbf$ is elliptically
symmetric [\citet{Eat86}], and is approximately satisfied when $p$
is large [\citet{HalLi93}]. Interestingly, as we show in
Section~\ref{sectionnonlinearha}, this assumption is no longer needed
for the unbiasedness in the more general setting of nonlinear
sufficient dimension reduction.

Here we note that, though Theorem \ref{theoremmain} is far from a
trivial generalization, the type of argument used in the proof is
somewhat standard in the SDR literature. See, for example,
\citet{LiDua89}, \citet{Coo98} and \citet{CooLi02}. It
is possible to extend the above theorem to more general objective
functions. For example, the theorem still holds if $a \mapsto a^+$ in
the objective function is replaced by any convex function $u(a)$.

\section{Estimation procedure for linear PSVM}
\label{sectionestimationprocedure}

\subsection{Estimation}\label{subsectionestimation}

We propose two ways to generate the set of pairs of slices
for PSVM. One, which we call ``left versus right'' (LVR),
repeatedly divides the predictors into two groups according
to a set of cutting points for the response. The other, which we call
``one versus another''
(OVA), partitions the predictors into several slices and pairs up all
possible slices.
We summarize the estimation procedure as follows.
\begin{longlist}[2$'$.]
\item[1.]
Compute the sample mean $\bar\Xbf$ and sample variance matrix $\hat
\Sigmabf$.
\item[2.] (LVR) Let $q_r, r=1, \ldots, h-1$, be $h-1$ dividing
points. For example, they can be equally spaced sample percentiles of
$\{Y_1, \ldots, Y_n \}$. Let
%
\begin{equation}\label{eq2slices}
\tilde Y_i^{r} = I(Y_i > q_r) - I(Y_i \le q_r)
\end{equation}
and let $(\hat{{\bolds\psi}}_r, \hat t_r )$, $r = 1, \ldots,
h-1$, be the
minimizer of
%
\begin{equation}\label{eqlvr}
{{\bolds\psi}}\trans\hat\Sigmabf{{\bolds\psi}}+
\lambda E_n \{1 - \tilde Y^{r}
[(\Xbf- \bar\Xbf)\trans{{\bolds\psi}}- t ]\}^+.
\end{equation}

\item[2$'$.] (OVA) Apply
SVM to each pair of slices from the $h$ slices. More
specifically, let $q_0 = \min\{Y_1, \ldots, Y_n \}$ and
$q_h = \max\{Y_1, \ldots, Y_n \}$. For each $(r,s)$
satisfying $1 \le r < s \le h$, let
\[
\tilde Y_i^{rs} =
I( q_{s-1} < Y_i \le q_s )- I( q_{r-1} < Y_i \le q_r).
\]
Let $(\hat{{\bolds\psi}}_{rs}, \hat t_{rs})$ be the minimizer
of the
objective function
\[
{{\bolds\psi}}\trans\hat\Sigmabf{{\bolds\psi}}+
\lambda E_n \{1 - \tilde Y^{rs}
[(\Xbf- \bar\Xbf)\trans{{\bolds\psi}}- t ]\}^+.
\]

\item[3.] Let $\hat{\mathbf v}_1, \ldots, \hat{\mathbf v}_d$ be the $d$ leading
eigenvectors of either one of the matrices
%
\begin{equation}\label{eqMn}
\hat\Mbf_n = \sum_{r=1}^{h-1} \hat{{\bolds\psi}}_r \hat
{{\bolds\psi}}{}\trans_r
\quad\mbox{or}\quad
\hat\Mbf_n = \sum_{r=1}^h \sum_{s = r+1}^h \hat{{\bolds\psi
}}_{rs} \hat{{\bolds\psi}}
{}\trans_{rs}.
\end{equation}
We use subspace
spanned by $\hat{\mathbf v}=( \hat{\mathbf v}_1, \ldots,
\hat{\mathbf v}_d)$ to
estimate $\spc_{Y|\Xbf}$.
\end{longlist}

Based on our experiences, LVR works best when the response is
a continuous variable, where $Y$ being larger or smaller has a
concrete physical meaning; OVA works best when the response
is categorical, where the values of $Y$ are simply labels of
classes, such as different vowels in our example in Section
\ref{sectiondataanalysis}.\vadjust{\goodbreak} Our numerical studies also
suggest that the estimation results are not overly sensitive
to the choice of the number of slices $h$, though a~larger~$h$ often
works better.

Standard packages for SVM minimize the objective function
(\ref{eqclue}) instead of (\ref{eqlvr}). However, they\vspace*{1pt} can
be modified to suit our procedure. Let $\zetabf=
\hat\Sigmabf{}\half{{\bolds \psi}}$ and $\Zbf =
\hat\Sigmabf{}\mhalf(\Xbf- \bar\Xbf)$. Then (\ref{eqlvr}) becomes
%
\begin{equation}\label{eqzetaproblem}
\zetabf\trans\zetabf+ \lambda E_n [1 - \tilde Y^r ( \Zbf\trans
\zetabf- t) ]^+.
\end{equation}
We can\vspace*{1pt} apply standard packages to minimize
(\ref{eqzetaproblem}) to obtain $\hat\zetabf$, whose transformation
$\hat\Sigmabf{}\mhalf\hat\zetabf$ is the desired minimizer of (\ref
{eqlvr}). We use the \texttt{kernlab} package in R to solve problem
(\ref{eqzetaproblem}). See \citet{KarMey06} for an exposition of
this package.

\subsection{Order determination}\label{subsectionorderdetermination}

Estimating the dimension $d$ of the central subspace is a vital
ingredient of sufficient dimension reduction estimation. Here, we
propose a cross-validated BIC procedure [\citet{Sch78}]
for this purpose. The BIC component of this procedure is an extension
of a criterion introduced
by \citet{WanYin08}, and is also related to \citet{ZhuMiaPen06}.
We refer to this combined procedure as CVBIC.

Let $\hat\Mbf_n$ be one of the matrices in (\ref{eqMn}), and
let $\lambda_i(\hat\Mbf_n)$ be its $i$th largest eigenvalue.
Let
$
G_n (k) = \sum_{i=1}^k \lambda_i ( \hat\Mbf_n) - c_1 (n) c_2 (k)$,
where $c_1(n)$ is a sequence of positive numbers or random
variables that converge\vspace*{2pt} (in probability) to 0, and $c_2(k)$ is
a nonrandom increasing function of $k$. Let
$\hat d$ be the maximizer of $G_n (k)$ over $\{0, \ldots, p\}$.
In Section \ref{sectionasymptoticanalysis}, we show that $P(\hat d =
d) \to1$.
The standard choices of $c_1 (n)$ and $c_2 (k)$ are $c_1 (n)
\propto n\mhalf\log(n) $ and $c_2 (k) = k$, so that the
penalty term is
$c_0 n\mhalf\log(n) k$,
where $c_0>0$ is a constant (or random variable)
of order $O (1)$ [or $O_\Psub(1)$]. Since the
eigenvalues $\lambda_i(\hat\Mbf_n)$ may differ for different
problems, it is sensible to make $c_0$ comparable to their
magnitude. One reasonable choice is to make $c_0$ proportional to
$\lambda_1 (\hat\Mbf_n)$,
leading to the following BIC-type criterion:
%
\begin{equation}\label{eqgamma0criterion}
\sum_{i=1}^k \lambda_i (\hat\Mbf_n) - a \lambda_1 (\hat\Mbf_n)
n\mhalf\log(n) k.
\end{equation}

We now turn to the choice of $a$. Though this choice
does not affect the consistency of $\hat d$, it does affect its
finite-sample performance. Moreover, from our experience
this choice is also sensitive to $p$, $d$, and the regression model.
For these reasons, it is important to have a systematic way of choosing
$a$. The SVM used
in our setting suggests naturally the cross-validation, because the
former provides a set of labels to validate.
We outline the CVBIC procedure as follows, using LVR as an illustration.

First, divide the data into a training set and a testing
set, denoted by
\[
\{(\acute\Xbf_1, \acute Y_1), \ldots, (\acute\Xbf_{n_1}, \acute
Y_{n_1})\},\qquad
\{(\grave\Xbf_1, \grave Y_1), \ldots, (\grave\Xbf_{n_2}, \grave
Y_{n_2})\}.
\]
Apply the PSVM to the training set with dividing points $q_1, \ldots,
q_{h-1}$ to obtain a set normal vectors $\acute{{\bolds\psi}}_1,
\ldots, \acute{{\bolds\psi}}_{h-1}$. Let $\acute\Mbf_{n_1} = \sum
_{i=1}^{h-1} \acute{{\bolds\psi}}_i \acute{{\bolds\psi}}{}\trans_i$.
Second, for a fixed $a$, maximize the criterion
(\ref{eqgamma0criterion}), with $\hat\Mbf_n$ replaced by
$\acute{\Mbf}_{n_1}$, to obtain an integer $k$. Let $\acute{\mathbf v}_1,
\ldots, \acute{\mathbf v}_k$ be the $k$ leading eigenvectors of
$\acute\Mbf_{n_1}$ and transform the testing predictors $\grave\Xbf_i$
to $\grave\Xbf_i^{(k)} = (\acute{\mathbf v}_1, \ldots, \acute {\mathbf
v}_k)\trans \grave\Xbf_i$, $i = 1, \ldots, n_2$. Third, let $\grave L_i
= I(\grave Y_i > q_r) - I(\grave Y_i \le q_r)$ be the\vspace*{1pt} true
label of $\grave Y_i$ in the testing set. Apply SVM to
$(\grave\Xbf_1^{(k)}, \grave L_1), \ldots, (\grave\Xbf _{n_2}^{(k)},
\grave L_{n_2})$ to predict $\grave L_1, \ldots, \grave L_{n_2}$.
Repeat this process for all dividing points and record the total number
of misclassifications. The optimal $a$ is the one that minimizes the
total number of misclassifications. Finally, substitute the optimal $a$
into (\ref{eqgamma0criterion}) and maximize it again using the full
data to estimate $d$. In Section \ref{subsectionestimated}, we
investigate the numerical performance of CVBIC under a variety of
combinations of~$p$,~$d$,~$n$ and regression models.

\subsection{Special features of linear PSVM}\label{subsectiontwofeatures}

As we conclude the exposition of the linear PSVM, we mention some
special features of this method. One is that it shares the
similar limitation with SIR when dealing with regression functions
that are symmetric about the origin. If the regression
function is
$ f( \| \Xbf\| ) $,
then all slices of the form $\{ \Xbf_i \dvtx Y_i \in S \}$ are roughly
concentric spheres in $\R^p$, which no hyperplane in $\R^p$ can
separate. However, as we shall see in
Sections \ref{sectionnonlinearha} and \ref{sectionsimulationstudies},
this is remedied by the kernel PSVM, because when mapped into higher
dimensional feature space the slices become linear again.

Another is that when dealing asymmetric regression functions, the
linear PSVM tends to work better than SIR for the following reason.
Recall that SIR is based, roughly, on the principal components of
the slice mean vectors of the form $E(\Xbf|Y \in S) - E(\Xbf)$,
where $S$ is an interval in $\Omega_Y$. This determines that
it downweights the slice means near the center of the data cloud,
where the Euclidean norm of $E(\Xbf|Y \in S) - E(\Xbf)$ is smaller.
However, it is well known that the regression function $E(Y|\Xbf)$
tends to be more accurately estimated near the center of the data
cloud [see, e.g., \citet{KutNacNet04},
Section 2.4]. In comparison, the linear PSVM relies on the normals of the
separating hyperplanes of the slices, which does not downweight the
data near the center.
As we will see from our simulation studies in Section
\ref{sectionsimulationstudies}, this brings substantial
improvement to the estimate.
We should point out, however, that there is an important exception. As
shown in \citet{Coo07} and
\citet{CooFor08}, under the assumption that $Y$ has a~finite
support and $\Xbf|Y$ has a conditional multivariate
normal distribution where $\var(\Xbf|Y)$ is independent of $Y$, SIR is
the maximum likelihood estimate of the
central subspace. In this case, no regular estimate can be more
efficient than SIR.
The mentioned advantage of linear PSVM applies mainly to the forward
regression setting where the
conditional distribution of $\Xbf|Y$ is typically non-Gaussian.

\section{Kernel PSVM for nonlinear dimension reduction}
\label{sectionnonlinearha}

In this section, we extend the PSVM to nonlinear
sufficient dimension reduction as defined by~(\ref{eqnonlinear}).
We first develop the objective function by
generalizing the linear PSVM objective function (\ref{eqmodified}),
and then establish the unbiasedness of the proposed nonlinear PSVM
estimator.

Before proceeding further, we note that the function $\phibf$ in relation
(\ref{eqnonlinear}) is not unique in the strict sense, but is
unique modulo injective transformations. Again, the situation is
parallel to linear sufficient dimension
reduction problem (\ref{relmodel}), where $\etabf\trans\Xbf$ is only
unique modulo injective linear transformations. Any injective linear
transformation of $\etabf\trans\Xbf$ is an equivalent linear
predictor, because it does not
change the linear subspace.
Likewise, for nonlinear SDR,
any injective transformation of $\phibf$ is an equivalent sufficient predictor,
because it does not change conditional independence (\ref{eqnonlinear}).

Let $\H$ be a Hilbert space of functions of $\Xbf$. In analogy
to the linear objective function (\ref{eqmodified}), consider
$\Lambda\dvtx \H\times\R\to\R^+$ defined by
%
\begin{equation}\label{eqvarianceform}
\Lambda( \psi, t) = \var[ \psi(\Xbf) ] + \lambda E \bigl[1 - \tilde Y \bigl(
\psi(\Xbf) - E
\psi(\Xbf) - t \bigr) \bigr]^+,
\end{equation}
where $\tilde Y$ is as defined in
(\ref{eqtildey}). To see that this is indeed a generalization
of (\ref{eqmodified}), consider the bilinear form from
$\H\times\H$ to $\R$ defined by
$
b(f_1, f_2) = \cov[ f_1 (\Xbf), f_2 (\Xbf) ].
$
Under the assumption that the mapping
%
\begin{equation}\label{eqembedding}
\H\to L_2 (P_{\mathbf X}), \qquad f \mapsto f
\end{equation}
is continuous, the bilinear form $b$ induces a bounded and self-adjoint
operator $\Sigma\dvtx\H\to\H$ such that $\langle f_1, \Sigma f_2
\rangle_{\H} = b(f_1,f_2)$, where $\langle\cdot, \cdot\rangle_\H$ is
the inner product in~$\H$. See, for example, \citet{Con90},
Theorem 2.2, and \citet{FukBacJor03}.
The objective function (\ref{eqvarianceform}) can now be
rewritten as
%
\begin{equation}\label{eqcommonform}
\Lambda(\psi, t) = \langle\psi, \Sigma\psi\rangle_\H+ \lambda E \bigl[1
- \tilde Y \bigl( \psi(\Xbf) - E \psi(\Xbf) - t \bigr) \bigr]^+.
\end{equation}
Thus, $\Lambda( \psi, t)$ is a generalization of $L({{\bolds
\psi}}, t)$
with the matrix $\Sigmabf$ replaced by the operator $\Sigma$,
the linear function ${{\bolds\psi}}\trans\Xbf$ replaced by an
arbitrary
function~$\psi$ in $\H$, and the inner product in $\R^p$ replaced
by the inner product in $\H$.
For the usual kernel SVM, the population-level objective function
is
\[
\langle\psi, \psi\rangle_\H+ \lambda E \bigl[1 - \tilde Y \bigl( \psi(\Xbf
) - E
\psi(\Xbf) - t \bigr) \bigr]^+.
\]
Comparing with (\ref{eqcommonform}), we see a parallel modification
to the linear case. The significance of this
modification is further discussed in Section \ref{sectionestimationkernelha}.

We now establish that, if $(\psi^*,t^*)$ is the minimizer of
$\Lambda(\psi, t)$, then $\psi^*$ is necessarily a function of
the sufficient predictor $\phibf(\Xbf)$ in the nonlinear
problem problem~(\ref{eqnonlinear}).
This is a generalization of the notion \textit{unbiasedness} in the
linear setting.
Our definition of unbiasedness (\ref{eqlinearunbiased})
in the linear sufficient dimension reduction setting is equivalent to
%
\begin{equation}\label{eqlinearunbiased1}
\mbox{$[\Tbf(F_0)]\trans\Xbf$ is a linear function of $\etabf
\trans
\Xbf$.}
\end{equation}
It is the statement (\ref{eqlinearunbiased1}) that is more
readily generalized to the nonlinear sufficient dimension
reduction setting: we simply require $\psi$ to be a function
of the sufficient predictor $\phibf(\Xbf)$ in (\ref{eqnonlinear}).
The following definition makes this notion rigorous. For a generic
random element $\Ubf$, let $\sigma\{\Ubf\}$ denote the
$\sigma$-field generated by $\Ubf$.
%
\begin{definition}
A function $\psi\in\H$ is unbiased for nonlinear sufficient dimension
reduction (\ref{eqnonlinear}) if it has a version that is
measurable $\sigma\{\phibf(\Xbf)\}$.
\end{definition}

The reason that we only require a version of $\psi$ to be measurable
$\sigma\{\phibf(\Xbf)\}$ is that the
$L_2$-metric ignores measure zero sets.
%
\begin{theorem}\label{theoremunbiasedness} Suppose the mapping (\ref
{eqembedding}) is continuous and:
\begin{longlist}[2.]
\item[1.] $\H$ is a dense subset of $L_2(P_{\Xbf})$,
\item[2.] $Y \indep\Xbf| \phibf(\Xbf)$.
\end{longlist}
If $(\psi^*, t^*)$ minimizes (\ref{eqcommonform})
among all $(\psi,t) \in\H\times\R$, then
$\psi^* (\Xbf)$ is unbiased.
\end{theorem}

Condition 1 is satisfied by some commonly used
reproducing kernel Hilbert spaces. For example, if $\mathcal G$ is
a reproducing kernel Hilbert space based on the Gaussian radial
basis, then the collection of functions $\{c + g\dvtx c \in\R,
g \in{\mathcal G} \}$ is dense in $L_2 (P_\Xbf)$. See \citet{FukBacJor09}.

It is important to note that in this more general setting we no longer
require any linearity assumption that resembles the one assumed in
Theorem~\ref{theoremmain}. In contrast, the kernel sliced inverse
regression developed by \citet{Wu08} and \citet{WuLiaMuk},
and functional sliced inverse regression by \citet{HsiRen09} all
require a version of the linearity condition to hold in the reproducing
kernel Hilbert space.

The notion of unbiasedness for sufficient dimension reduction
is more akin to Fisher consistency than to unbiasedness
in the classical setting. While unbiasedness in the classical setting can
exclude many useful statistics, Fisher consistency often guarantees
correct asymptotic behavior without putting undue
restrictions on the expectation. Moreover, an estimator that is not
Fisher consistent is clearly undesirable,
because it is guaranteed \textit{not} to converge to the true parameter.
For these reasons unbiasedness for linear SDR
is a useful criterion, even though some useless estimators (such as
$\bblfnull
$) are unbiased.
Unbiasedness for nonlinear SDR plays the parallel role, except
that it only requires the estimator to be an arbitrary, rather than a
linear, function of the true predictor.
This relaxation also allows us to establish the unbiasedness of PSVM
without evoking the linearity condition.

Theorem \ref{theoremunbiasedness} assumes that $\Lambda(\psi,t)$
attains its minimum in $\H\times\R$. We think this is a reasonable
assumption for the following reasons. As shown below,~$\Lambda(\psi,t)$
is lower semicontinuous with respect to the weak topology in
$\H\times\R$. Since any closed, bounded, and convex set in a Hilbert
space is compact with respect to the weak topology [\citet{Wei80},
Theorem~4.25, \citet{Con90},  Corollary~V.1.5], by the generalized
Weierstrass theorem [\citet{KurZab05}, Section 7.3],
$\Lambda(\psi, t)$ attains its minimum within any such set in
$\H\times\R$. The next proposition establishes this fact. Let $\hxstar$
be the Hilbert space $\hx\times\R$ endowed with the inner product $
\langle\psi_1, \psi_2 \rangle_\hx+ t_1 t_2. $
%
\begin{proposition} If $\hx$ is an RKHS with its kernel
$\kappa$ satisfying $E   \kappa(\Xbf,\allowbreak \Xbf) < \infty$, then
$\Lambda(\psi,t)$ is lower semicontinuous with respect to the weak
topology in~$\H'$,
and attains its minimum in any closed, bounded, and convex
set in $\hxstar$.
\end{proposition}

\section{Estimation of kernel PSVM and invariant kernel}
\label{sectionestimationkernelha}

The purpose of this section is twofold. First, because we
have modified $\langle\psi, \psi\rangle_\H$ to $\langle
\psi, \Sigma\psi\rangle_\H$ in the kernel SVM objective
function, we can no longer use the standard SVM packages
to solve for $\psi^*$. Therefore, we reformulate the
minimization of~$\Lambda(\psi, t)$ as quadratic
programming that can be solved by available
computer packages. Second, in deriving this quadratic
programming problem, we gain more insights into the meaning
and significance of this modification. As we shall see, by
replacing $\langle\psi, \psi\rangle_\H$ by $\langle\psi,
\Sigma\psi\rangle_\H$, we are in effect making SVM
invariant with respect to the marginal distribution of $\Xbf$.
Intuitively, since we are using SVM to make inference about
the conditional distribution of~$Y|\Xbf$, it is plausible
that the procedure does not depend on the marginal
distribution of $\Xbf$.

Let $\H$ be a linear space of functions from
$\Omega_\Xbf$ to $\R$ spanned by
$
{\mathcal F}_n = \{ \psi_1, \ldots, \psi_k \}.
$
The choice of these functions will be
discussed later, but it will ensure
$E_n [ \psi_i (\Xbf) ] = 0$, so that $\psi_i ({\mathbf x}) =
\psi_i ({\mathbf x}) -
E_n \psi_i (\Xbf)$.
Let
%
\begin{equation}\label{eqcappsi}
\Psibf=
\pmatrix{
\psi_1 (\Xbf_1) & \cdots& \psi_k (\Xbf_1) \cr
\vdots& \ddots& \vdots\cr
\psi_1 (\Xbf_n)& \cdots&\psi_k(\Xbf_n)}.
\end{equation}
Then the sample
version of the objective function (\ref{eqcommonform}) is
%
\begin{equation}\label{eqhatlambda}
\hat\Lambda(\cbf) = n^{-1} \cbf\trans\Psibf\trans\Psibf\cbf+
\lambda n^{-1} \sum_{i=1}^n [ 1 - \tilde Y_i ( \Psibf_i\trans\cbf-
t) ]^+,
\end{equation}
where $\Psibf_i\trans= ( \psi_1 (\Xbf_i), \ldots, \psi_k (\Xbf
_i) )$
and $\cbf\in\R^k$.
We minimize $\hat\Lambda(\cbf)$ among all~$\cbf$.

In the following, $\tilde\ybf=(\tilde y_1, \ldots, \tilde y_n)\trans$
and $\alphabf,
\betabf, \xibf\in\R^n$. The symbol
$\le$ represents componentwise inequality. The symbol $\odot$
represents the Hadamard
product between matrices. For a matrix $\Abf$ of full
column rank, $\Pbf_{\Abf}$ is the projection $\Abf
(\Abf\trans\Abf)^{-1} \Abf\trans$. The symbols
$\zero$ and $\one$ represent, respectively, the
$n$-dimensional vectors whose entries are 0 and 1.
%
\begin{theorem}\label{theoremquadraticprogramming}
$\!\!\!$If $\cbf^*$ minimizes $\hat\Lambda(\cbf)$ over $\R^k$,
then $\cbf^* = \frac{1}{2} (\Psibf\trans\Psibf)^{-1}\Psibf\trans
(\tilde\ybf\odot\alphabf^*)$, where $\alphabf^*$ is the
solution to the quadratic programming problem:
%
\begin{eqnarray}\label{eqtargetproblem}
&&\mbox{maximize } \one\trans\alphabf- \tfrac{1}{4}( \alphabf
\odot\tilde\ybf)\trans\Pbf_{ \Psibf} (\alphabf\odot\tilde\ybf)
\nonumber\\[-8pt]\\[-8pt]
&&\qquad\mbox{subject to } \zero\le\alphabf\le\lambda\one,
\alphabf\trans\tilde\ybf= 0.
\nonumber
\end{eqnarray}
\end{theorem}

Note that the quadratic programming problem
(\ref{eqtargetproblem}) differs from that of the
standard kernel SVM, where the projection
$\Pbf_{ \Psibf}$ is replaced by the kernel
matrix $\Kbf_n = \{\kappa(i,j)\dvtx i, j = 1, \ldots, n \}$
for some positive definite bivariate mapping $\kappa\dvtx
\Omega_\Xbf\times\Omega_\Xbf\to\R$. The kernel matrix
$\Kbf_n$ uniquely determines the sample estimate of the
covariance operator $\Sigma$, which bears the information
about the shape of the marginal distribution of $\Xbf$.
By replacing $\Kbf_n$ with~$\Pbf_\Psibf$, we are, in
effect, removing the information about $\Xbf$. For this reason we call
the matrix $\Pbf_{ \Psibf}$ an \textit{invariant kernel}.

For the function class ${\mathcal H}$, we use the reproducing
kernel Hilbert space based on the mapping
$\kappa$. Common
choices of $\kappa$ include the polynomial kernel $\kappa
({\mathbf x}_1, {\mathbf x}
_2) =
({\mathbf x}_1 \trans{\mathbf x}_2 + c)^r$, where $r$ is a
positive integer, and the
Gaussian radial kernel $\kappa({\mathbf x}_1,\allowbreak {\mathbf x}_2) =
e^{-\gamma\|{\mathbf x}_1 -
{\mathbf x}_2\|^2}$, where $\gamma> 0$.
Let
%
\begin{equation}\label{eqrkhs}
\H_\kappa=\{c_0 + c_1 \kappa(\cdot, \Xbf_1) + \cdots+ c_n \kappa
(\cdot, \Xbf_n)\dvtx c_0, \ldots, c_n \in\R\}
\end{equation}
with inner product specified by
$\langle\kappa(\cdot, {\mathbf a}), \kappa(\cdot, {\mathbf b})
\rangle= \kappa({\mathbf a}, {\mathbf b})$. In the standard
kernel SVM, it is a common practice to use all functions in $\H_\kappa$
as $\H$. However, the invariant
nature of our kernel, $\Pbf_\Psibf$, determines that we
cannot use all those functions, because if so then
$\Pbf_\Psibf$ becomes nearly an identity matrix (note that if $\Pbf
_{{\bolds\psi}}$ were an
identity matrix then the objective function in (\ref{eqtargetproblem})
would become independent of $\Xbf_1, \ldots, \Xbf_n$).
We instead use the principal functions of the linear
operator $\Sigma_n$, as defined by
$
\langle\psi_1, \Sigma_n \psi_2 \rangle= \cov_n [ \psi_1 (\Xbf),
\psi
_2 (\Xbf) ]
$, as our basis ${\mathcal F}_n$. Here $\cov_n(\cdot, \cdot)$ denotes
sample covariance.
Let $\Qbf_n = \Ibf_n - \Jbf_n/n$, where $\Ibf_n$ is the $n\times n$
identity matrix and $\Jbf_n$ is
the $n \times n$ matrix whose entries are 1.
The next proposition tells us how to find the eigenfunctions
of $\Sigma_n$. Its proof is easy and omitted.
%
\begin{proposition}\label{propositioneigenfunction}
Let $\wbf= (w_1, \ldots, w_n)$, $\psi_{\wbf}= \sum w_i [\kappa
({\mathbf x},
\Xbf_i) - E_n \kappa({\mathbf x}, \Xbf)]$.
The following statements are equivalent:
\begin{enumerate}
\item$\wbf$ is an eigenvector of the matrix $\Qbf_n \Kbf_n \Qbf_n$
with eigenvalue $\lambda$;
\item$\psi_{\wbf}$ is an eigenfunction of the operator $\Sigma_n$ with
eigenvalue $\lambda/n$.
\end{enumerate}
If $\lambda\ne0$, then either statement implies $(\psi_{\wbf} (\Xbf
_1), \ldots, \psi_{\wbf} (\Xbf_n)) = \lambda\wbf\trans$.
\end{proposition}

Although the eigenvectors of $\Qbf_n \Kbf_n \Qbf_n$ and the eigenfunctions
of $\Sigma_n$ are similar objects, it is the latter that can be
evaluated at any ${\mathbf x}$, not just the observed $\Xbf_1,
\ldots, \Xbf_n$. This property is important for prediction.
Essentially, we use the first $k$ eigenfunctions $\phi_1, \ldots,
\phi
_k$ of $\Sigma_n$ as
the functions in ${\mathcal F}_n$. This is equivalent to using
$\{a_1 \phi_1, \ldots, a_k \phi_k\} \equiv\{\psi_1, \ldots, \psi
_k\}$
for any nonzero
$a_1, \ldots, a_k$. We choose $a_i$ to satisfy
$a_i(\psi_i(\Xbf_1), \ldots, \psi_i(\Xbf_n))\trans= \wbf_i$, where
$\wbf_i$
is the eigenvector of $\Qbf_n \Kbf_n \Qbf_n$, corresponding to its
$i$th eigenvalue $\lambda_i$.
Thus $a_i = 1/\lambda_i$. With this choice, $\Psibf$ is simply $(\wbf
_1, \ldots, \wbf_k)$.
The choice of number of basis functions, $k$, should allow
sufficient flexibility but not as large as $n$; our
experiences indicate that the choice of $k$ in the range
$n/3 \sim2n/3$ works well.
We summarize the kernel PSVM estimation procedure as follows.

\begin{longlist}[1.]
\item[1.] (Optional) Marginally standardize $\Xbf_1, \ldots, \Xbf_n$.
Let $\hat\mu_r$ and $\hat\sigma_r^2$ be the sample mean and
sample variance $X_{1r}, \ldots, X_{nr}$. Reset $X_{ir}$ to be
$(X_{ir} - \hat\mu_r) / \hat\sigma_r$. The purpose of this step
is so that the kernel $\kappa$ treats different components of
$\Xbf_i$ more or less equally. This step can be omitted if the
components of $\Xbf_i$ have similar variances.

\item[2.] Choose a kernel $\kappa$ and the number of basis functions
$k$ (say $k=n/2$). Compute $\Psibf= (\wbf_1, \ldots, \wbf_k)$ and
$\Pbf
_\Psibf$ from $\Qbf_n \Kbf_n \Qbf_n$.
\item[3.] Divide the sample according to LVR or OVA, each yielding
a set of slices. For each pair of slices, solve the quadratic
programming problem in Theorem \ref{theoremquadraticprogramming}
using the $\Pbf_\Psibf$ computed from step 2. This gives\vspace*{-1pt}
coefficient vectors~$\cbf_1^*, \ldots, \cbf_{\tilde h}^* \in\R^k$,
where $\tilde{h} = h-1$ for LVR and $\tilde{h} = {h \choose2}$ for
OVA.
\item[4.] Compute the first $d$ eigenvectors, ${\mathbf v}_1, \ldots,
{\mathbf v}_d$,
of the matrix $\sum_{s=1}^{\tilde{h}} \cbf_s^* {\cbf_s^*}\trans$.
Denote the $r$th component of of ${\mathbf v}_s$ as $v_{sr}$.
\item[5.] The $s$th sufficient predictor evaluated at ${\mathbf x}$ is
$v_{s1} \psi_1 ({\mathbf x}) + \cdots+ v_{sk} \psi_k
({\mathbf x})$, where
$\psi_r ({\mathbf x}) = \lambda_r^{-1}\sum_{i=1}^n w_{ri} [
\kappa
({\mathbf x}, \Xbf_i) - E_n \kappa({\mathbf x}, \Xbf)]$.
If step 1 is used, then ${\mathbf x}$ should be marginally
standardized by the $\hat\mu_r$ and $\hat\sigma_r$ computed
from that step.
\end{longlist}
Many computing packages are available to solve the quadratic programming
problem in step 3. We use the \texttt{ipop} program in the
\texttt{kernlab} package in R. See \citet{Karetal04}.
If the Gaussian radial kernel is used in step~2, then
we recommend choosing $\gamma$ as
%
\begin{equation}\label{eqsamplegamma}
\gamma= 1/\tau^2, \qquad \tau= \frac{1}{{n \choose2}}\sum_{i < j,
j=2}^n\| \Xbf_i - \Xbf_j \|.
\end{equation}
Alternatively, we can use the population version of the above quantity,
%
\begin{equation}\label{eqpopulationgamma}
\gamma= 1/ ( E \| \Xbf- \Xbf' \| )^2,
\end{equation}
where $\Xbf$ and $\Xbf'$ are independent $N(\zero,\Ibf_p)$ random vectors.
This quantity can be easily evaluated by Monte Carlo.
In Section \ref{sectionsimulationstudies}, we use (\ref
{eqpopulationgamma})
for large-scale simulations to avoid repeated evaluations of
(\ref{eqsamplegamma}), whereas in Section~\ref{sectiondataanalysis}
we use (\ref{eqsamplegamma}) for the real data analysis, where it needs
to be calculated only once.
Some authors recommend sample median
in (\ref{eqsamplegamma}). See
\citet{Greetal05} and
\citet{FukBacJor09}.
This does not make a~significant difference in our examples.

\section{Asymptotic analysis of linear PSVM}
\label{sectionasymptoticanalysis}

In this section, we give a comprehensive asymptotic analysis of linear
PSVM estimator introduced in Sections \ref{sectionlinearha} and
\ref{sectionestimationprocedure}. This is developed in three parts.
First, we derive the influence function for the normal vector
$\hat{{\bolds\psi}}$ based on two slices. In this part, we employ some
asymptotic properties of SVM developed recently by
\citet{JiaZhaCai08}. In the second part, we derive the asymptotic
distribution of the linear PSVM estimator, $(\hat{\mathbf v}_1, \ldots,
\hat{\mathbf v}_d)$, defined in Section \ref {subsectionestimation}. In
the third part, we establish the consistency of the order determination
criterion introduced in Section \ref{subsectionorderdetermination}.

\subsection{Influence function for support vector machine}

The asymptotic results of \citet{JiaZhaCai08} are largely
applicable here except for three places: our SVM involves an
additional $\Sigmabf$; our $\lambda$ is fixed but the $\lambda$
in their paper depends on $n$; they did not derive the explicit
form of the hessian matrix---and hence neither the asymptotic
variance---but we are interested in the explicit asymptotic
distribution. The first two points are minor but the third needs
nontrivial additional work. We only consider the case where~$\tilde Y$
is defined through a partition $\{A_1, A_2\}$ of
$\Omega_Y$. Thus, our results only apply to the LVR scheme. The
asymptotic analysis the OVA scheme can be carried out similarly,
and is omitted.

We first develop some notation. Let\vspace*{1pt} $\thetabf=
({{\bolds\psi}}\trans, t)\trans$, $\Zbf= (\Xbf\trans, \tilde Y)\trans$,
$\Xbf^* = ( \Xbf\trans$, $-1)\trans$, $\Sigmabf^* = \diag(\Sigmabf , 0)$.
Then
%
\begin{equation}\label{eqmfunction}
{{\bolds\psi}}\trans\Sigmabf{{\bolds\psi}}+ \lambda
[1 - \tilde Y ( \Xbf\trans{{\bolds\psi}}- t )]^+ = \thetabf
\trans\Sigmabf
^* \thetabf- \lambda( 1-
\thetabf\trans\Xbf^* \tilde Y )^+.
\end{equation}
We denote this function by $m(\thetabf, \Zbf)$. Let $\Omega_\Zbf$ be
the support of $\Zbf$ and let
${\mathbf h}\dvtx \Theta\times\Omega_\Zbf\to\R^r$
be a function of $(\thetabf, \Zbf)$. Let $D_\thetabf$
denote the $(p+1)$-dimensional column vector of differential operators
$(\partial/ \partial\theta_1, \ldots, \partial/ \partial
\theta_{p+1})\trans$.
The
next theorem gives the gradient of the support vector machine
objective function $E[ m (\thetabf, \Zbf)]$.
%
\begin{theorem}\label{theoremgradient} Suppose,
for each $ \tilde y = -1, 1$, the distribution of
$\Xbf| \tilde Y = \tilde y $ is dominated by the Lebesgue measure
and $E(\| \Xbf\|^2) < \infty$.
Then
%
\begin{equation}\label{eqfirst-order}
D_\thetabf E [ m (\thetabf, \Zbf)] = (2 {{\bolds\psi}}\trans
\Sigmabf, 0)\trans
-\lambda E [ \Xbf^* \tilde Y I ( 1 - \thetabf\trans\Xbf^*\tilde Y >
0 )].
\end{equation}
\end{theorem}

We now present the hessian matrix of support vector
machine, which leads to the asymptotic variance of $\hat\thetabf$.
To our knowledge, this is the first time that the asymptotic
variance is explicitly given. This result is then used to derive
the asymptotic distribution of the linear PSVM estimator.
%
\begin{theorem}\label{theoremhessian} Suppose $\Xbf$ has a convex
and open support and its conditional distributions given $\tilde Y =1$
and $\tilde Y = -1$ are dominated by the Lebesgue measure. Suppose,
moreover:\vfill\eject
\begin{longlist}[1.]
\item[1.] for any linearly independent
${{\bolds\psi}}, \deltabf\in\R^p$, $\tilde y = -1, 1$, and
$v \in\R$,
the following function is continuous:
\[
u \mapsto E ( \Xbf^* | {{\bolds\psi}}\trans\Xbf= u, \deltabf
\trans\Xbf= v,
\tilde Y = \tilde y ) f_{{{\bolds\psi}}\trans
\Xbf| \deltabf\trans\Xbf, \tilde Y} (u | v, \tilde y );
\]
\item[2.] for any $i\,{=}\,1, \ldots, p$, and $ \tilde y\,{=}\,-1, 1$, there is a
nonnegative function~$c_i(v, \tilde y)$ with
$E[ c_i(V, \tilde Y) | \tilde Y ] < \infty$ such that
\[
v E( X_i | {{\bolds\psi}}\trans\Xbf= u, \deltabf\trans\Xbf
= v, \tilde Y =
\tilde y )
f_{{{\bolds\psi}}\trans\Xbf| \deltabf\trans
\Xbf, \tilde Y} (u
| v, \tilde y ) \le c_i (v, \tilde y);
\]
\item[3.] there is a nonnegative function
$c_0 (v, \tilde y)$ with $E [c_0 (V,\tilde Y) | \tilde Y ] < \infty$
such that $ f_{{{\bolds\psi}}\trans\Xbf|
\deltabf\trans\Xbf,
\tilde Y} (u | v, \tilde y ) \le c_0 (v, \tilde y)$.
\end{longlist}
Then the function $\thetabf\mapsto D_\thetabf E[ m(\thetabf, \Zbf)]$
is differentiable in all directions with derivative matrix
\[
2 \diag(\Sigmabf,0) + \lambda\sum_{\tilde y=-1,1}P(\tilde Y =
\tilde y )
f_{{{\bolds\psi}}\trans\Xbf|\tilde Y}
(t + \tilde y | \tilde y)
E (\Xbf^* {\Xbf^*}\trans| {{\bolds\psi}}\trans\Xbf= t +
\tilde y ).
\]
Furthermore, if the function \mbox{$({{\bolds\psi}}, t)\,{\mapsto}\,
f_{{{\bolds\psi}}\trans\Xbf|\tilde Y}(t\,{+}\,\tilde y |\tilde y)
E (\Xbf^* {\Xbf^*}\trans| {{\bolds\psi}}\trans\Xbf\,{=}\,t\,{+}\,\tilde y )$} is
continuous, then $D_\thetabf[m(\thetabf, \Zbf)]$ is jointly
differentiable with respect to $\thetabf$.
\end{theorem}

Joint differentiability and directional differentiability are sometimes
refered to as Frechet differentiability and Gateaux differentiability.
The latter
is generally weaker than the former. In a finite-dimensional space, having
continuous directional derivative in all directions implies joint
differentiability
[\citet{Bicetal93}, page 453].
The next theorem gives the influence function for support vector
machine.
%
\begin{theorem}\label{theoreminfluence}
If the conditions in Theorems \ref{theoremgradient} and \ref{theoremhessian}
are satisfied, then
\[
\hat\thetabf= \thetabf_0 - \Hbf\inv\{ (2 {{\bolds\psi
}}_0\trans\Sigmabf, 0
)\trans- \lambda E_n [ \Xbf^* \tilde Y
I ( 1 - \tilde Y \thetabf_0\trans\Xbf^* > 0 ) ]\} + o_{\Psub}
(n\mhalf),
\]
where $\Hbf$ is hessian matrix given by Theorem \ref{theoremhessian}.
\end{theorem}

The proof is similar to that of \citet{JiaZhaCai08} and
is omitted. Alternatively, one can prove it by applying Theorem 5.23 of
\citet{van98}.

\subsection{\texorpdfstring{Asymptotic distribution of $(\hat\vb_1, \ldots, \hat \vb_d)$}
{Asymptotic distribution of (v_1,...,v_d)}}

Consider\vspace*{1pt} a fixed division point $q_r$, where $r \in\{1, \ldots, h-1\}$.
Let $\tilde Y^r$ be as defined in (\ref{eq2slices}),
and $\Zbf^r=(\Xbf\trans, \tilde Y^r)\trans$. Let $\thetabf_{0 r}
= ({{\bolds\psi}}_{0 r}\trans, t_{0 r})\trans$ be the
minimizer of
$E [m(\thetabf, \Zbf^r)]$,
and $\hat\thetabf_{r} =
(\hat{{\bolds\psi}}{}\trans_{r}, \hat t_{r})\trans$ be the
minimizer of
$E_n [ m (\thetabf, \Zbf^r)]$.
Let $\Hbf_r$ be the hessian matrix of
$E[ m (\thetabf, \Zbf^r)]$, and let $\Fbf_r$ be the first $p$ rows
of $\Hbf_r^{-1}$. By Theorem \ref{theoreminfluence},
%
\begin{equation}\label{eqinfluenceforsvmir}
\hat{{\bolds\psi}}_r = {{\bolds\psi}}_{0 r} - \sbf_r
(\thetabf_{0 r}, \Zbf^r) +
o_{\Psub} (n\mhalf),
\end{equation}
where
$
\sbf_r (\thetabf, \Zbf^r) = \Fbf_r [(2 {{\bolds\psi}}\trans
\Sigmabf, 0 )\trans
- \lambda\Xbf^* \tilde Y^r
I ( 1 - \tilde Y^r \thetabf\trans\Xbf^* > 0 ) ].
$
Let
\[
\hat\Mbf_n = \sum_{r=1}^{h-1} \hat{{\bolds\psi}}_r \hat
{{\bolds\psi}}{}\trans_r,\qquad
\Mbf_0 = \sum_{r=1}^{h-1} {{\bolds\psi}}_{0 r} {{\bolds
\psi}}_{0 r}\trans.
\]
For a matrix $\Abf\in\R^{r_1 \times r_2}$, let $\Kbf_{r_1, r_2}
\in\R^{r_1 r_2 \times r_1 r_2}$ be the commutation matrix defined by
the relation $\Kbf_{r_1, r_2} \vec(\Abf) = \vec(\Abf\trans)$. See
\citet{MagNeu79}. Two properties of $\Kbf_{r_1,r_2}$ that will
prove useful for our purpose are that $\Kbf_{r_1,r_2} = \Kbf_{r_2,
r_1}\trans$ and that for any $\Bbf\in\R^{r_3 \times r_4}$,
%
\begin{equation}\label{eqcommutationrule}
\Abf\otimes\Bbf= \Kbf_{r_1, r_3} ( \Bbf\otimes\Abf) \Kbf_{r_4, r_2}.
\end{equation}
We now present the asymptotic distribution of $\hat\Mbf_n$.
%
\begin{theorem}\label{theoremcandidatematrix} Under the assumptions in
Theorems \ref{theoremgradient} and \ref{theoremhessian}, $\sqrt n
\vec(\hat\Mbf_n - \Mbf_0 )$
converges to multivariate normal with mean $\zero$ and variance
\[
(\Ibf_{p^2} + \Kbf_{p,p}) \sum_{r=1}^{h-1} \sum_{t=1}^{h-1}
( {{\bolds\psi}}_{0 r}{{\bolds\psi}}_{0 t}\trans
\otimes\Lambdabf_{rt} ) (\Ibf_{p^2} +
\Kbf_{p,p}),
\]
where
$
\Lambdabf_{rt} = E [\sbf_r (\thetabf_{0 r}, \Zbf^r)\sbf_t\trans
(\thetabf_{0 t}, \Zbf^t)]$.
\end{theorem}

This result leads directly to the asymptotic distribution of $\hat\Vbf=
(\hat{\mathbf v}_1, \ldots,\allowbreak \hat{\mathbf v}_d )$. Since, by
Theorem \ref{theoremmain},
$\sp(\Mbf_0 )
\subseteq\spc_{Y|\Xbf}$, we have $\rank(\Mbf_0) \le d$.
We make the working assumption that
$\rank(\Mbf_0) = d$. This means we exclude the situations
where the regression surface is symmetric about the origin.
Since $\Mbf_0$ is positive semi-definite, it has the
spectral decomposition $\Ubf\Dbf\Ubf\trans$, where~$\Ubf$
is a $p \times d$ matrix whose columns are the eigenvectors of $\Mbf_0$
corresponding to nonzero eigenvalues, and $\Dbf$ is a
$d \times d$ diagonal matrix with diagonal elements being
the nonzero eigenvalues. The following corollary
is a direct consequence of Theorem \ref{theoremcandidatematrix} and
\citet{BurPfe08}.
Its proof is omitted.

\begin{corollary}\label{corollaryvhat} Under the assumptions in Theorems
\ref{theoremgradient} and \ref{theoremhessian}
and\break $\operatorname{rank} (\Mbf_0) = d$,
$
\sqrt n   \vec( \hat\Vbf- \Vbf_0 ) \cid N (
\bblfnull,
\Upsilonbf),
$
where $\Upsilonbf$ is the $pd \times pd$ matrix
\[
(\Dbf\inv\Ubf\trans\otimes\Ibf_p)(\Ibf_{p^2} + \Kbf_{p,p}) \sum
_{r=1}^{h-1} \sum_{t=1}^{h-1}
( {{\bolds\psi}}_{0 r}{{\bolds\psi}}_{0 t}\trans
\otimes\Lambdabf_{rt} ) (\Ibf_{p^2} +
\Kbf_{p,p})
( \Ubf\Dbf\inv\otimes\Ibf_p).
\]
\end{corollary}

It is possible to refine the PSVM estimator by introducing weights to~$\hat\Mbf_n$.
Take the LVR scheme for example. Let $\Psibf= (\hat
{{\bolds\psi}}_1, \ldots, \hat{{\bolds\psi}}_{h-1})$. Let $\Abf$ be an
${h-1}$ by ${h-1}$ matrix. Rather than working with $\hat\Mbf_n$, we
could base the spectral decomposition on a weighted matrix $\hat\Mbf_n
(\Abf) = \Psibf\trans\Abf\Psibf$. Let $\hat{\mathbf v}(\Abf) =
(\hat{\mathbf v}_1 (\Abf), \ldots, \hat{\mathbf v}_d (\Abf))$ be the
first $d$ eigenvectors of $\hat\Mbf_n (\Abf)$. One way to determine the
optimal $\Abf$ is by minimizing a real-valued monotone function (say
trace) of the asymptotic variance matrix of $\vec[\hat{\mathbf
v}(\Abf)]$, which can be extracted from the asymptotic distribution.
This type of argument was used in Li (\citeyear{Li00},
\citeyear{Li01}) to construct optimal estimating equations.
Alternatively, one can develop an optimal procedure using the minimum
distance approach introduced by \citet{CooNi05}. We leave these to
future research.

\subsection{Consistency of the BIC-type criterion}\label
{subsectionconsistencyofd}

In the following, we say a~sequence of random variables $W_n$
converges in probability to infinity (\mbox{$W_n \cip\infty$})
if, for any $K > 0$, $\lim_{n\to\infty} P(|W_n| > K) = 1$. Let $\hat d$
is the maximizer of~$G_n (k)$ over $\{0, \ldots, p\}$ as defined in Section
\ref{subsectionorderdetermination}.
%
\begin{theorem}\label{theorembicconsistency} Suppose $P(c_1 (n) > 0)=
1$, $ c_1 (n) \cip0$,
$ n\half c_1 (n) \cip\infty$,\break and~$c_2 (k)$ is an increasing
function of $k$. Under the conditions in
Theorems~\ref{theoremmain},~\ref{theoremgradient} and~%
\ref{theoremhessian} and $\rank(\Mbf_0) = d$, we have $\lim_{n \to
\infty} P ( \hat d = d) = 1$.
\end{theorem}

Note that we have again made the working assumption
$\rank(\Mbf_0) = d$. However, even when this assumption is
violated the theorem still holds with~$d$ replaced by the
rank of $\Mbf_0$.

\section{Simulation studies}\label{sectionsimulationstudies}

In this section, we compare the linear and kernel PSVM with four other
methods based on the idea of inverse regression: SIR, the sliced
average variance estimator [SAVE; \citet{CoWe91}], directional
regression [DR; \citet{LiWan07}], and kernel sliced inverse regression
[\citet{Wu08}]. We also investigate the performance of the CVBIC for order
determination.

\subsection{Linear dimension reduction}\label{subsectionlineardr}

We use the following models:
\begin{eqnarray*}
\mbox{Model I:}\quad Y&=&{X_1}/{[0.5+(X_2 +1)^2]} + \sigma\varepsilon,
\\
\mbox{Model II:}\quad Y&=& X_1 (X_1+X_2+1) +\sigma\varepsilon, \\
\mbox{Model III:}\quad Y &=& ( X_1^2 + X_2^2)^{1/2}     \log( { X_1^2 +
X_2^2} )^{1/2} + \sigma\varepsilon,
\end{eqnarray*}
where $\Xbf\sim N(0, \Ibf_p)$, $p = 10, 20, 30$, $\varepsilon\sim
N(0, 1)$ and $\sigma= 0.2$. The sample size~$n$ is taken to be $100$.
The first two models, which are taken from \citet{Li91N2}, are asymmetric
about 0, but the last one is symmetric about 0. As we have discussed
in Section \ref{subsectiontwofeatures}, linear PSVM, like SIR, does
not work when the regression surface is symmetric about 0. The first
two examples show how the linear PSVM compares with other methods in
the situations where it works. The purpose of the last model is to
provide a benchmark of error when it fails, so that we can gauge how
the kernel PSVM improves the situation in the next comparison.\vadjust{\goodbreak}

To evaluate the performance of each method, we use the distance measure
suggested by \citet{LiZhaChi05}. Specifically, let
$\mathcal{S}_1$ and $\mathcal{S}_2$ be two subspaces of $\R^p$. Then
%
\begin{equation}\label{eqlzc}
\dist(\spc_1, \spc_2) = \| \Pbf_{\spc_1} - \Pbf_{\spc_2} \|,
\end{equation}
where $\Pbf_{\mathcal{S}_1}$ and $\Pbf_{\mathcal{S}_2}$ are orthogonal
projections on to $\mathcal{S}_1$ and $\mathcal{S}_2$, and
\mbox{$\|\cdot\|$} is a matrix norm. In the following the Frobenius
norm is used.

For SAVE and DR, we use $h=4$ slices, and for SIR, we use $h=8$ slices,
having roughly the same number of points. Our choices of $h$ are in
line with the usual practice in the SDR literature for such a sample
size. For methods such as SAVE and DR that involve the second-order
inverse moment, $h$ is suggested to be chosen smaller than that for
methods such as SIR which only involve the first-order inverse moment
[\citet{LiZhu07}]. For the linear PSVM, the cost $\lambda$ is
taken to be 1. The number of division points ($q_r$) is 20. We have
tried some other numbers of division points and obtained very similar
results. In general, our experiences suggest that a relatively large
number of dividing points is preferable. The results are presented in
Table \ref{table1}. The entries are of the form $a (b)$ where $a$ is
the mean, and $b$ is the standard deviation of the distance criterion
(\ref{eqlzc}) calculated from 200 simulated samples. The last row in
Table \ref{table1} records the CPU time (in seconds) each method uses
for Model I with $p = 10$ (on a Dell OptiPlex 745 desktop computer with
speed 2.66 GHz).

\begin{table}
\caption{Estimated means and simulation standard errors (in
parentheses) of the distance
measure (\protect\ref{eqlzc}) and mean computation times (in second) of linear
sufficient dimension reduction methods}\label{table1}
\begin{tabular*}{\tablewidth}{@{\extracolsep{\fill}}l r c c c c@{}}
\hline
\textbf{Models} & $\bolds{p}$ & \textbf{SIR} & \textbf{SAVE} &
\textbf{DR} & \textbf{Linear PSVM} \\
\hline
\hphantom{II}I & $10$ & 0.84 (0.22)& 1.55 (0.19)& 1.02 (0.23)& 0.65 (0.17) \\
& $20$ & 1.14 (0.18)& 1.93 (0.05)& 1.32 (0.17)& 0.93 (0.16) \\
& $30$ & 1.31 (0.14)& 1.96 (0.03)& 1.48 (0.11)& 1.17 (0.14) \\
[4pt]
\hphantom{I}II & $10$ & 1.20 (0.27)& 1.43 (0.16)& 1.17 (0.23)& 0.85 (0.25) \\
& $20$ & 1.51 (0.19)& 1.72 (0.15)& 1.46 (0.14)& 1.26 (0.23) \\
& $30$ & 1.67 (0.16)& 1.84 (0.12)& 1.63 (0.12)& 1.58 (0.17) \\
[4pt]
III & $10$ & 1.80 (0.13)& 0.87 (0.21)& 0.85 (0.20)& 1.65 (0.16) \\
& $20$ & 1.89 (0.08)& 1.46 (0.20)& 1.45 (0.20)& 1.85 (0.10) \\
& $30$ & 1.93 (0.05)& 1.72 (0.12)& 1.71 (0.12)& 1.93 (0.05) \\
[4pt]
Time & & 0.03\hphantom{(00\ .0)} & 0.01\hphantom{(00\ .0)}
& 0.03\hphantom{(00\ .0)} & 0.16\hphantom{(00\ .0)} \\
\hline
\end{tabular*}
\end{table}

Table \ref{table1} shows that the linear PSVM consistently performs
better than the other methods in all cases for models I and II. The
intuition behind this improvement is explained in Section
\ref{subsectiontwofeatures}. Also, as expected, the linear PSVM and SIR
do not perform well for Model III because of the symmetry of the
regression function. However, as we will see in the next comparison,
this defect is no longer present in the kernel PSVM. The linear PSVM
requires more computing time than the classical methods, mainly because
it needs to process more dividing points, and, for each dividing point,
the full data (rather than a slice of data) are processed.

\subsection{Nonlinear dimension reduction}

As we have mentioned, Model III is symmetric about 0, and the linear
PSVM fails. To a certain degree, the shape of regression surface of
Model II is also symmetric about 0. We now use these two models to
investigate the performance of the kernel PSVM for nonlinear sufficient
dimension reduction. In terms of linear dimension reduction, Model III
has two sufficient predictors, $X_1$, $X_2$, but in terms of nonlinear
dimension reduction, it has only one sufficient predictor, $(X_1^2 +
X_2^2)^{1/2}$, or any monotone function of it. The kerenl PSVM is
designed to recover a~monotone transformation of $(X_1^2 +
X_2^2)^{1/2}$ without having to assume any regression model. In doing
so, it solves two problems at one stroke---further reducing the
dimension from 2 to 1, and avoiding the difficulty of SIR in dealing
with symmetric responses.

\begin{figure}

\includegraphics{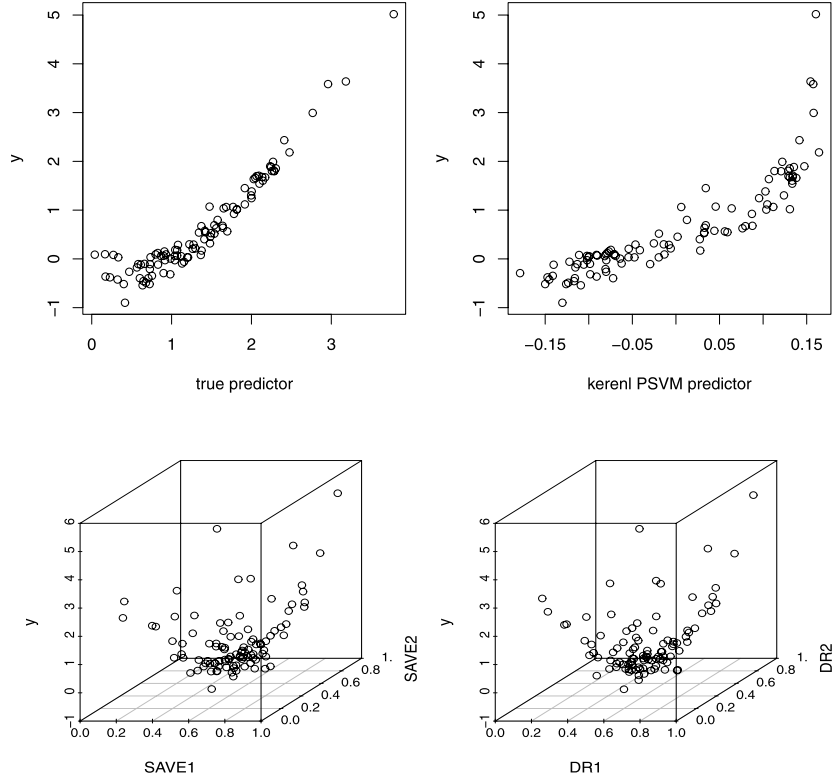}

\caption{Comparison between linear and nonlinear sufficient dimension
reduction methods. Upper left panel: true nonlinear predictor
$\sqrt{X_1^2 + X_2^2}$ versus $Y$; upper right panel: first (nonlinear)
PSVM predictor versus $Y$; lower left: first two SAVE predictors versus~$Y$;
lower right panel: first two DR predictors versus~$Y$.}
\label{figuremodel8}
\vspace*{-3pt}
\end{figure}

To illustrate the idea, in Figure \ref{figuremodel8} we present the
2-D and 3-D scatter plots for $Y$ versus the nonlinear and linear
predictors obtained by different methods. The upper left panel is the
2-D scatter plot for $Y$ versus the true nonlinear predictor
$(X_1^2 + X_2^2)^{1/2}$; the upper right panel is the 2-D scatter
plot of $Y$ versus the first kerenl PSVM predictor; the
lower panels are \mbox{3-D} scatter plots for $Y$ versus the first two
predictors from SAVE and DR. We can see that all three methods
capture the right shape of the regression function, but kernel PSVM only
requires one predictor and its sufficient plot appears sharper
(bearing in mind that the upper right panel only has to resemble a
\textit{monotone transformation} of the upper left panel).

To make a more precise comparison, we need to design a new criterion
that can compare one nonlinear predictor with two linear predictors;
the criterion (\ref{eqlzc}) is no longer suitable for this purpose.
Since the nonlinear sufficient predictor estimates a monotone
function of $(X_1^2 + X_2^2)^{1/2}$, we use the absolute value of
Spearman's correlation to measure their closeness, which is invariant
under monotone transformation [\citet{KutNacNet04},
page 87]. To measure the closeness between two linear predictors and
the true nonlinear predictor $(X_1^2 + X_2^2)^{1/2}$, let
$
(U_{11}, \ldots, U_{1n}), (U_{21}, \ldots, U_{2n})
$
represent the two linear predictors obtained by SAVE or DR. These
predictors estimate linear combinations of $X_{1i}, X_{2i}$ but do
not specify $X_{1i}$, $X_{2i}$ themselves. We therefore regress
$T_i = X_{1i}^2 + X_{2i}^2$ on
\[
\{(1, U_{1i}, U_{2i}, U_{1i}^2, U_{2i}^2)\dvtx i = 1, \ldots, n \}.
\]
If $U_{1i}$ and $U_{2i}$
are (linearly independent) linear combinations of $X_{1i}$ and~$X_{2i}$,
then this regression is guaranteed to recover the true predictor
$T_i$ regardless of the specific form of
the linear combinations. Let $\hat T_1, \ldots, \hat T_n$ be the fitted
responses of this regression. We use the absolute values of Spearman's
correlation between $T_i$ and $\hat T_i$ to measure the performance
of SAVE and DR.

\begin{table}
\caption{Estimated means and simulation standard errors (in
parentheses) of Spearman correlations
of linear and nonlinear sufficient dimension reduction}\label{table2}
\begin{tabular*}{\tablewidth}{@{\extracolsep{\fill}}lcccccccc@{}}
\hline
& \multicolumn{4}{c}{\textbf{Model II}} &
\multicolumn{4}{c@{}}{\textbf{Model III}} \\ [-4pt]
& \multicolumn{4}{c}{\hrulefill} &
\multicolumn{4}{c@{}}{\hrulefill} \\
$\bolds{p}$ & \textbf{SAVE} & \textbf{DR} & \textbf{KSIR}
& \textbf{KPSVM} & \textbf{SAVE} & \textbf{DR} & \textbf{KSIR}
& \textbf{KPSVM}\\
\hline
10 &0.53& 0.67& 0.88& 0.92& 0.79 & 0.79 & 0.89 & 0.90
\\
&(0.13)& (0.11)& (0.07)& (0.02)& (0.09) & (0.08) & (0.05) & (0.02) \\
[4pt]
20 &0.37& 0.53& 0.68& 0.86& 0.56 & 0.57 & 0.59 & 0.81
\\
&(0.14)& (0.09)& (0.17)& (0.03)& (0.11) & (0.11) & (0.18) & (0.03) \\
[4pt]
30 &0.30& 0.43& 0.55& 0.83& 0.47 & 0.48 & 0.42 & 0.77
\\
&(0.13)& (0.10)& (0.23)& (0.04)& (0.11) & (0.11) & (0.21) & (0.04) \\
\hline
\end{tabular*}
\end{table}

We compute these numbers for 200 simulation samples, and tabulate
their means and standard deviations in Table \ref{table2}.
Note that large numbers represent better performance, and all
numbers are between 0 and 1. The SAVE and DR estimators are
computed in exactly the same way as in the linear dimension
reduction comparison. For the kernel PSVM, the cost is 1, the number
of division points is still~20, the kernel is the Gaussian radial
basis, and the number of principal eigenfunctions of $\Sigma_n$ is
taken to be 60. The parameter $\gamma$ is calculated by
(\ref{eqpopulationgamma}), which are approximately $0.0526$,
$0.0257$ and $0.0169$ for $p=10,20,30$, respectively. We see that
the kernel PSVM actually performs better than SAVE and DR, even
though it uses only one predictor. It also performs better than KSIR.
Moreover, the accuracy of the kernel PSVM remains reasonably high for
larger $p$, where the accuracies of SAVE,
DR, and KSIR drop considerably.\vadjust{\goodbreak}

\subsection{Estimation of structural dimension}
\label{subsectionestimated}

We now investigate the performance of the CVBIC order-determination procedure
for a variety of combinations of $(p,d,n)$.
We still use models I and II, but,
to include different $d$, we add the following models which both have
$d =1$:
\begin{eqnarray*}
\mbox{Model IV:}\quad Y&=&{X_1}/{[0.5+(X_1 +1)^2]} + \sigma\varepsilon,\\
\mbox{Model V:}\quad Y&=&X_1 (2 X_1 +1) +\sigma\varepsilon.
\end{eqnarray*}
These are derived from models I and II by replacing
$X_2$ by $X_1$.

We apply CVBIC in conjunction with PSVM to Models I, II, IV and V, with
$(d,n,p)$ ranging over the set $\{1,2\} \times\{200, 300, 400, 500 \}
\times\{ 10, 20, 30 \}$. The training and testing sample sizes are $n_1
= n_2 = n/2$. We take 20 dividing points $q_r$ as equally-spaced sample
quantiles of $\acute Y_1, \ldots, \acute Y_{n_1}$. As a comparison we
also apply the order-determination procedure for SIR based on Theorem~%
5.1 of \citet{Li91N2} with significant level $\alpha= 0.05$. The
results are presented in Table \ref{table3}, where the entries are the
percentage of correct estimation of $d$ out of 200
%
\begin{table}
\caption{Rate of correct order determination by SIR and PSVM in $\%$}\label{table3}
\begin{tabular*}{\tablewidth}{@{\extracolsep{\fill}}lcccd{3.0}cd{3.0}
cd{3.0}cc@{}}
\hline
&&&
\multicolumn{2}{c}{$\bolds{n=200}$} & \multicolumn{2}{c}{$\bolds{n=300}$}&
\multicolumn{2}{c}{$\bolds{n=400}$}& \multicolumn{2}{c@{}}{$\bolds{n=500}$} \\[-4pt]
&&&
\multicolumn{2}{r}{\hrulefill} & \multicolumn{2}{r}{\hrulefill}
& \multicolumn{2}{r}{\hrulefill}
& \multicolumn{2}{c@{}}{\hrulefill}\\
\textbf{Model} & \multicolumn{1}{c}{$\bolds{d}$}
& \multicolumn{1}{c}{$\bolds{p}$} & \multicolumn{1}{c}{\textbf{SIR}}
& \multicolumn{1}{c}{\textbf{PSVM}} & \multicolumn{1}{c}{\textbf{SIR}}
& \multicolumn{1}{c}{\textbf{PSVM}} & \multicolumn{1}{c}{\textbf{SIR}}
& \multicolumn{1}{c}{\textbf{PSVM}} & \multicolumn{1}{c}{\textbf{SIR}}
& \multicolumn{1}{c@{}}{\textbf{PSVM}}\\
\hline
\hphantom{II}I & 1 & 10 & 92 & 96 &92 & 100& 97 &
100 &98 & 100 \\
& & 20 & 80 & 82 &95 & 96& 96 & 100 &94 & 100 \\
& & 30 & 65 & 54 &92 & 94& 94 & 98 &96 & 100 \\
[4pt]
\hphantom{I}II & 2 & 10 & 67 & 80 &86 & 85& 97 & 90
&98 & \hphantom{0}94 \\
& & 20 & 36 & 64 &66 & 84& 85 & 86 &96 & \hphantom{0}84 \\
& & 30 & 22 & 32 &55 & 77& 73 & 84 &88 & \hphantom{0}80 \\
[4pt]
IV & 1 & 10 & 39 & 82 &60 & 84& 71 & 93
&75 & \hphantom{0}95 \\
& & 20 & 28 & 78 &42 & 74& 54 & 76 &68 & \hphantom{0}80 \\
& & 30 & 15 & 78 &33 & 84& 44 & 78 &60 & \hphantom{0}80 \\
[4pt]
\hphantom{I}V & 2 & 10 & 93 &100 &96 & 100& 96 &
100 &97 & 100 \\
& & 20 & 94 & 98 &96 & 100& 96 & 100 &96 & 100 \\
& & 30 & 96 & 98 &96 & 99& 96 & 100 &92 & 100 \\
\hline
\end{tabular*}
\end{table}
simulated samples for each of the 48 combinations of $(\mathrm{model},
p, n )$. Table \ref{table3} shows that CVBIC works very well, with
percentage of correct estimation reaching as high as 100\% for sample
size of 200 (training sample size 100). In almost all cases, PSVM
compares favorably with SIR for order determination. Also clear from
the table is the trend of increasing accuracy for both methods as $n$
increases.

\section{Application and further discussions}\label{sectiondataanalysis}

$\!\!\!$We now compare the kernel PSVM with SIR, SAVE, and DR in a real data analysis
concerning recognition of vowels. The data can be found in the
UCI Machine Learning Repository
(\url{http://archive.ics.uci.edu/ml/datasets}). The response variable
$Y$ is a categorical variable of 11~levels, representing
different vowel sounds. The predictor $\Xbf$ is a 10-dimensional
vector describing the features of a sound.
For clear presentation, we select only three vowels: the sounds
in \textit{heed}, \textit{head} and \textit{hud}, with training and testing
sample sizes being 144 and 126, respectively.

\begin{figure}

\includegraphics{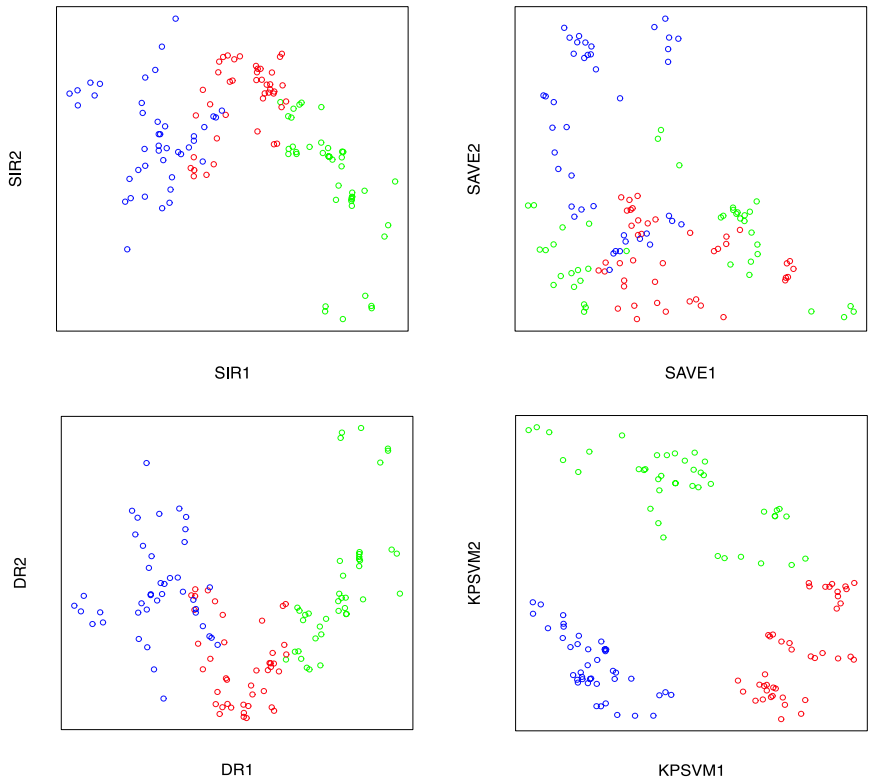}

\caption{First two predictors based on SIR, SAVE, DR, and the kernel PSVM
plotted for the vowel recognition testing data set. Green, red and blue
colors indicate the vowel sounds in \textup{heed}, \textup{head}, \textup{hud}.}
\label{figurevow2d}
\end{figure}

For each dimension reduction method, we find a set of
sufficient predictors from the training data, and evaluate them
at the testing set, resulting in a~sufficient plot for the testing set.
Given that the testing data are independent of the training data
from which the sufficient predictors is derived, the degree of separation
of the vowels in the sufficient plot objectively reflects the
discriminating power of a dimension reduction method.
The four scatter plots in Figure \ref{figurevow2d} present
the first two predictors found by SIR (upper left panel), SAVE
(upper right panel), DR (lower left panel), and the kerenl PSVM
(lower right panel). For the kernel PSVM, the OVA scheme is used.
The basis functions are the first 40 eigenfunctions of the
operator $\Sigma_n$ derived from the Gaussian radial kernel,
whose parameter $\gamma$ is calculated by~(\ref{eqsamplegamma}).
The cost $\lambda$ is 1. We have varied the
number of eigenfunctions (from~10 to 60) and the cost (from 0.5
to 20), but they do not seem to result in significant difference
in the degree of separation in the test data.

From Figure \ref{figurevow2d}, we see that the kernel PSVM achieves
much better separation of the three vowels in the test data than
the other three methods. The second best performer is DR, followed
by SIR and SAVE. It is also interesting to note that the various
degrees of separation are also reflected in the sufficient plots;
that is, the distance between \textit{heed} and \textit{hud} is larger
than those between \textit{heed} and \textit{head}, and \textit{head} and
\textit{hud}.

We would like to comment that classification,
though important, is not the sole purpose for sufficient dimension
reduction, and that linear and nonlinear
sufficient dimension reductions
have their own strengths in reducing, discriminating, visualizing, and
interpreting high-dimensional data.
To illuminate the point, consider an example where variation, rather
than location,
is the differentiating characteristic.
Let $Y$ be a bernoulli variable with $P(Y=1) = P(Y=0) = 1/2$
and
\[
(\Xbf| Y =y) \sim N
\left( \bblfnull,
\pmatrix{
\sigma^2 (y) \Ibf_2 & 0 \cr
0 & \Ibf_{p-2}}\right),
\]
where $\sigma^2 (0) = 1$ and $\sigma^2(1) = 10$.
Let $(\Xbf_1, Y_1), \ldots, (\Xbf_n, Y_n)$ be a sample from this model,
where $n = 200$ and $p = 10$.
For simplicity, we fix the number of cases of $Y=1$ at
$n/2$, because this has no bearing on our problem.
In this case, the central subspace is $\sp(\ebf_1, \ebf_2)$,
where $\ebf_i = (0,\ldots,1,\ldots,0)\trans$ with the $1$ occupying the
$i$th position.

%
\begin{figure}

\includegraphics{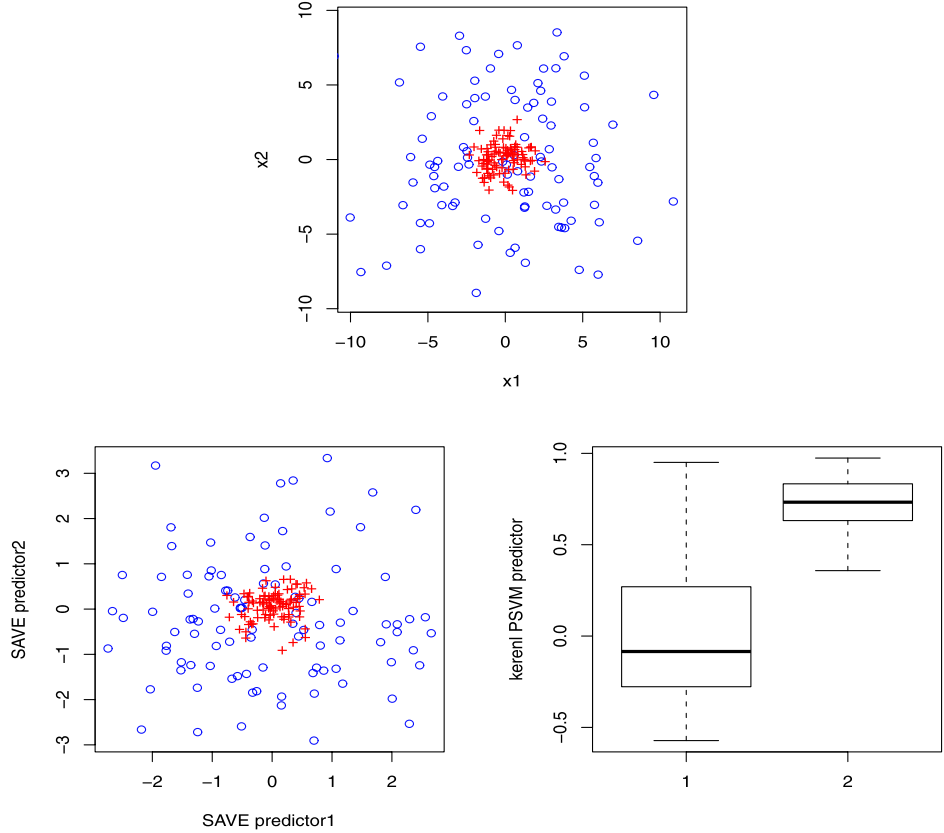}

\caption{Variation as the differentiating characteristic.
Blue $\circ$ represents the $Y=0$ cases and red $+$ represents the
$Y=1$ cases.}\label{figureaaa}
\end{figure}

We apply SAVE and the kernel PSVM and the results are shown in Figure~%
\ref{figureaaa}, where the top panel shows the scatter plot for
the true sufficient predictors $X_1$ and~$X_2$, the lower left panel
shows the first two SAVE predictors, and the lower right panel shows
the boxplot of a single kernel PSVM predictor.
Since for a single variable
we cannot produce a scatter plot, for clarity we use a boxplot to
represent the predictor. The value of the kernel PSVM predictor
is represented
by the height in the boxplot; the two boxes represents the two groups.
All three plots are based on the testing data.
What is interesting is that kernel PSVM in some sense ``translates''
the difference in variation into the
difference in location. The intuitive reason is that there is a
quadratic---and hence variance---component
in the kernel mapping, but in the mapped high-dimensional space the
variance component is treated as an augmented part of
feature vector [as in $(x,x^2)$]. Of course this is only a
simplification of the situation:
there is still significant difference
in variation in the kernel PSVM predictor for the two groups.

In this case, linear dimension reduction methods such as SAVE have a~definite advantage,
both for their clear separation of variation and for their good
interpretability.
In the meantime, this example also shows that kernel PSVM is capable of
differentiating variation, to the degree
comparable to SAVE, but its interpretability is not as direct as SAVE.

\begin{figure}

\includegraphics{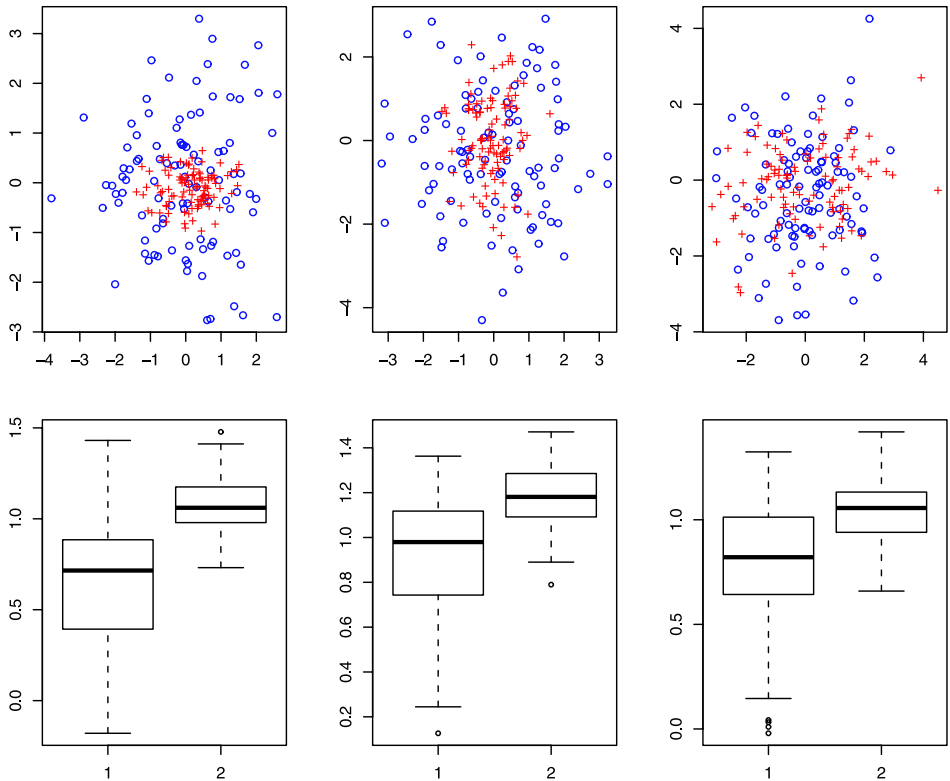}

\caption{Degrees of separation by SAVE (upper panels) and kernel PSVM
(lower panels) for higher dimensions: $p=60$ (left panels),
$p = 80$ (middle panels) and $p = 100$ (right panels).}
\label{figurereferee1comment103}
\vspace*{-3pt}
\end{figure}

Another desirable feature of the kernel PSVM is that
its accuracy is more stable than the classical methods as the dimension
$p$ increases. Figure \ref{figurereferee1comment103}
shows
the sufficient predictors derived from SAVE and kernel PSVM for $p =
60, 80, 100$ (from left to right). The upper panels
are the scatter plots for the first two SAVE predictors, and the lower
panels are the boxplots representing
the single kernel PSVM predictor.
Again, all
plots are based on testing data. We see that SAVE gradually loses its
discriminating power as $p$ is increased
to 100, whereas the discriminating power of kernel PSVM remains
reasonably strong.

\section*{Acknowledgments}

We are very grateful to three referees and an Associate Editor, whose
many useful comments and suggestions greatly broadened and deepened
an earlier version of this work.\vadjust{\goodbreak}



\printaddresses

\end{document}